\documentclass[11pt,a4paper]{article}
  \def\SvZfontcode{8}
  \def\SvZslantedGreekCapitals{1}
  \usepackage{amsmath,amsthm,amscd}

\def\SvZrequireslantedRedef{0}
\ifcase\SvZfontcode
\usepackage{bm}
\usepackage{amssymb}
\def\SvZrequireslantedRedef{1}
\or
\usepackage{fouriernc}
\usepackage{bm}
\def\SvZrequireslantedRedef{1}
\or
\usepackage{fourier}
\usepackage{bm}
\def\SvZrequireslantedRedef{1}
\or
\usepackage{amssymb}
\if\SvZslantedGreekCapitals 1
\usepackage[slantedGreek]{mathptmx}
\else
\usepackage{mathptmx}
\fi
\DeclareMathAlphabet{\bm}{OT1}{ptm}{b}{it} 
\or
\usepackage{kmath,kerkis}
\usepackage{bm}
\def\SvZrequireslantedRedef{1}
\or
\if\SvZslantedGreekCapitals 1
\usepackage[sc,slantedGreek]{mathpazo}
\else
\usepackage[sc]{mathpazo}
\fi
\usepackage{bm}
\or
\usepackage[adobe-utopia]{mathdesign}
\usepackage{bm}
\or
\usepackage[urw-garamond]{mathdesign}
\usepackage{bm}
\or
\usepackage[bitstream-charter]{mathdesign}
\usepackage{bm}
\or
\usepackage{lmodern}
\usepackage{bm}
\def\SvZrequireslantedRedef{1}
\or
\usepackage{arev}
\fi

\if\SvZrequireslantedRedef 1
\if\SvZslantedGreekCapitals 1

\renewcommand{\Gamma}{\varGamma}
\renewcommand{\Delta}{\varDelta}
\renewcommand{\Theta}{\varTheta}
\renewcommand{\Lambda}{\varLambda}
\renewcommand{\Xi}{\varXi}
\renewcommand{\Pi}{\varPi}
\renewcommand{\Sigma}{\varSigma}
\renewcommand{\Upsilon}{\varUpsilon}
\renewcommand{\Phi}{\varPhi}
\renewcommand{\Psi}{\varPsi}
\renewcommand{\Omega}{\varOmega}
\fi
\fi

\renewcommand{\phi}{\varphi}

\reversemarginpar

\usepackage{url}
\usepackage{listings}
\usepackage{kbordermatrix}
\usepackage{booktabs}
\usepackage{graphicx}
\usepackage{color}

\ifx\SvZinPresentation\undefined
\usepackage{float}
\usepackage{longtable}
\usepackage[margin=10pt,labelfont=bf,textfont=it]{caption}

\usepackage[colorlinks]{hyperref}

\renewcommand{\theenumi}{\textit{\roman{enumi}}}
\renewcommand{\labelenumi}{(\theenumi)}

\newcommand{\mathds}{\mathbb}

\DeclareMathOperator{\rank}{rank}

\newcommand{\Q}{ \mathds{Q} }
\newcommand{\R}{ \mathds{R} }

\newcommand{\Z}{ \mathds{Z} }
\newcommand{\C}{ \mathds{C} }
\newcommand{\N}{\mathds{N}}
\newcommand{\field}{\mathds{F}}
\newcommand{\group}{G}

\DeclareMathOperator{\GF}{GF}
\newcommand{\ring}{R}
\newcommand{\parf}{\ensuremath{\mathbb{P}}} 
\newcommand{\reg}{\uniform_0}
\newcommand{\dyadic}{\mathds{D}}
\newcommand{\nreg}{\uniform_1}
\newcommand{\psru}{\mathds{S}}
\newcommand{\sru}{\ensuremath{\sqrt[6]{1}}}
\newcommand{\golrat}{\mathds{G}}
\newcommand{\splittable}{\mathds{Y}}
\newcommand{\ger}{\mathds{GE}}
\newcommand{\gauss}{\hydra_2}
\newcommand{\hydra}{\mathds{H}}
\newcommand{\cyclo}{\mathds{K}}

\newcommand{\uniform}{\mathds{U}}

\DeclareMathOperator{\cra}{cr}
\DeclareMathOperator{\crat}{Cr}

\DeclareMathOperator{\assoc}{Asc}
\DeclareMathOperator{\fun}{{\cal F}}
\newcommand{\lift}{\mathds{L}}

\renewcommand{\tilde}{\widetilde}
\renewcommand{\hat}{\widehat}

\newcommand{\ignore}[1]{}

\DeclareMathOperator{\sign}{sgn}

\let\Oldsetminus\setminus

\newcommand{\delete}{\ensuremath{\!\Oldsetminus\!}}
\newcommand{\contract}{\ensuremath{\!/}}

\newcommand{\minorof}{\ensuremath{\preceq}}

\newcommand{\bip}{G}

\ifx\SvZinPresentation\undefined
\newtheorem{theorem}{Theorem}[section]

\newtheorem{lemma}[theorem]{Lemma}
\newtheorem{proposition}[theorem]{Proposition}
\newtheorem{definition}[theorem]{Definition}
\newtheorem{corollary}[theorem]{Corollary}

\newtheorem{claim}{Claim}[theorem] 
\newtheorem{aclaim}{Claim}[theorem] 
\newtheorem{sclaim}{Claim}[claim] 

\newenvironment{claimenv}{\list{}{\rightmargin0pt\leftmargin10pt\topsep0pt}\item[]}{\endlist}

\newenvironment{subproof}{\begin{claimenv}\begin{proof}}{\end{proof}\end{claimenv}}
\fi
\newtheorem{conjecture}[theorem]{Conjecture}
\newtheorem{question}[theorem]{Question}
 
\begin{document}
  \title{Lifts of matroid representations over partial fields}
  \author{R. A. Pendavingh and S. H. M. van Zwam\thanks{E-mail: \url{rudi@win.tue.nl}, \url{svzwam@win.tue.nl}. This research was supported by NWO, grant 613.000.561.}}
  \maketitle
  \center{Dedicated to Lex Schrijver on the occasion of his sixtieth birthday.}
  \abstract{
    There exist several theorems which state that when a matroid is representable over distinct fields $\field_1, \ldots, \field_k$, it is also representable over other fields. We prove a theorem, the Lift Theorem, that implies many of these results.

    First, parts of Whittle's characterization of representations of ternary matroids follow from our theorem. Second, we prove the following theorem by Vertigan: if a matroid is representable over both $\GF(4)$ and $\GF(5)$, then it is representable over the real numbers by a matrix such that the absolute value of the determinant of every nonsingular square submatrix is a power of the golden ratio. Third, we give a characterization of the 3-connected matroids having at least two inequivalent representations over $\GF(5)$. We show that these are representable over the complex numbers.

    Additionally we provide an algebraic construction that, for any set of fields $\field_1, \ldots, \field_k$, gives the best possible result that can be proven using the Lift Theorem.
  }

\section{Introduction}
Questions regarding the representability of matroids pervade matroid theory. They underly some of the most celebrated results of the field, as well as some tantalizing conjectures. A famous theorem is the characterization of regular matroids due to Tutte. We say that a matrix over the real numbers is \emph{totally unimodular} if the determinant of every square submatrix is in the set $\{-1,0,1\}$.
\begin{theorem}[Tutte \cite{Tut65}]\label{thm:reg}
Let $M$ be a matroid. The following are equivalent:
\begin{enumerate}
  \item $M$ is representable over both $\GF(2)$ and $\GF(3)$;
  \item  $M$ is representable by a totally unimodular matrix;
  \item $M$ is representable over every field.
\end{enumerate}
\end{theorem}
Whittle~\cite{Whi95,Whi97} proved very interesting results of a similar nature. Here is one example. We say that a matrix over the real numbers is \emph{totally dyadic} if the determinant of every square submatrix is in the set $\{0\}\cup\{\pm 2^k \mid k \in \Z\}$.
\begin{theorem}[Whittle~\cite{Whi97}]\label{thm:dyadicintro}
Let $M$ be a matroid. The following are equivalent:
\begin{enumerate}
  \item $M$ is representable over both $\GF(3)$ and $GF(5)$;
  \item $M$ is representable by a totally dyadic matrix;
  \item $M$ is representable over every field that does not have characteristic 2.
\end{enumerate}
\end{theorem}
A third example is the following result. We say that a matrix over the real numbers is \emph{golden ratio} if the determinant of every square submatrix is in the set $\{0\}\cup\{\pm \tau^k \mid k \in \Z\}$. Here $\tau$ is the golden ratio, i.e. the positive root of $x^2-x-1 = 0$.
\begin{theorem}[{Vertigan}]\label{thm:golratintro}
Let $M$ be a matroid. The following are equivalent:
\begin{enumerate}
  \item $M$ is representable over both $\GF(4)$ and $GF(5)$;
  \item $M$ is representable by a golden ratio matrix;
  \item  $M$ is representable over $\GF(p)$ for all primes $p$ such that $p = 5$ or $p \equiv \pm 1 \mod 5$, and also over $\GF(p^2)$ for all primes $p$.
\end{enumerate}
\end{theorem}
The common feature of these theorems is that representability over a set of finite fields is characterized by the existence of a representation matrix over some field such that the determinants of square submatrices are restricted to a certain set $S$.
Semple and Whittle~\cite{SW96} generalized this idea. They introduced \emph{partial fields}: algebraic structures where multiplication is as usual, but addition is not always defined. The condition ``all determinants of square submatrices are in a set $S$'' then becomes ``all determinants of square submatrices are defined''. In this paper we present a general theorem on partial fields from which results like Theorems~\ref{thm:reg}--\ref{thm:golratintro} follow. We employ a mixture of combinatorial and algebraic techniques.

We start our paper, in Section~\ref{sec:pfintro}, with a summary of the work of Semple and Whittle~\cite{SW96}. We note here that we have changed the definition of what it means for a sum to be \emph{defined}, because with the definition proposed by Semple and Whittle a basic proposition, on which much of their work is based, is false.
We give numerous additional definitions and basic results, and introduce notation to facilitate reasoning about representation matrices of a matroid.
The ideas behind our definitions are ubiquitous --- they capture the way Truemper~\cite{Tru92} relates matroids and representation matrices, they occur in Section 6.4 of Oxley~\cite{oxley}, and even the ``representative matrices associated with a dendroid'' in Tutte~\cite{Tuthom} are essentially the same thing. 

Section~\ref{sec:lifts} contains the main theorem of this paper, the Lift Theorem (Theorem~\ref{thm:lift}). It gives a sufficient condition under which a matroid that is representable over a partial field $\parf$ is also representable over a partial field $\hat\parf$. The condition is such that it can be checked for classes of matroids as well. 

In Section~\ref{sec:examples} we give applications of the Lift Theorem. First we give alternative proofs for a significant part of Whittle's~\cite{Whi97} characterization of the ternary matroids that are representable over some field of characteristic other than 3. We also prove Vertigan's Theorem~\ref{thm:golratintro} and two new results, namely a characterization of the 3-connected matroids that have at least two inequivalent representations over $\GF(5)$, and a characterization of the subset of these that is also representable over $\GF(4)$.

Another result by Vertigan, Theorem~\ref{thm:pfinring}, states that every partial field can be seen as a subgroup of the group of units of a commutative ring. We give a proof of this theorem in Section~\ref{sec:rings}.
We show that a matroid representable over some partial field is in fact representable over a field. This complements the theorem by Rado~\cite{Ra57} that every matroid representable over a field is also representable over a finite field. We also show that for every partial field homomorphism there exists a ring homomorphism between the corresponding rings.

We use these insights to define a ring and corresponding partial field for which, by construction, the premises of the Lift Theorem hold. With this partial field we can formulate a result like Theorems~\ref{thm:reg}--\ref{thm:golratintro} for any finite set of finite fields. We show that our construction gives the ``best possible'' partial field to which the Lift Theorem applies.

Finally we present, in Section~\ref{sec:questions}, a number of unsolved problems that arose during our investigations.

In a related paper~\cite{PZ08conf} we show that in some instances the Lift Theorem can be pushed a little further. In particular we show that for a 3-connected matroid $M$ it may happen that only a sub-partial field is needed to represent $M$.

The statements of Theorems~\ref{thm:golratintro} and \ref{thm:pfinring} were mentioned in  Geelen et al.~\cite{GOVW98} and in Whittle~\cite{Whi05} as unpublished results of Vertigan. This work was started because we wanted to understand Vertigan's results. Our proofs were found independently. Vertigan informs us that he had, in fact, proven Lemma~\ref{lem:liftpf}, using methods very similar to those found in Section~\ref{sec:lifts} of this paper, and that he had deduced Theorem~\ref{thm:golratintro} from that.

\section{Preliminaries}\label{sec:pfintro}
\subsection{Notation}
If $S$, $T$ are sets, and $f:S\rightarrow T$ is a function, then we define
\begin{align}
  f(S):= \{f(s) \mid s \in S\}.
\end{align}
We denote the restriction of $f$ to $S'\subseteq S$ by $f|_{S'}$. We may simply write $e$ instead of the singleton set $\{e\}$.

If $S$ is a subset of elements of some group, then $\langle S \rangle$ is the subgroup generated by $S$. If $S$ is a subset of elements of a ring, then $\langle S \rangle$ denotes the \emph{multiplicative} subgroup generated by $S$. All rings are commutative with identity. The group of elements with a multiplicative inverse (the \emph{units}) of a ring $\ring$ is denoted by $\ring^*$. As usual, if $S$ is a set of indeterminates, then $\ring[S]$ denotes the polynomial ring over $\ring$.

Our graph-theoretic notation is mostly standard. All graphs encountered are simple. We use the term \emph{cycle} for a simple, closed path in a graph, reserving \emph{circuit} for a minimal dependent set in a matroid. An undirected edge (directed edge) between vertices $u$ and $v$ is denoted by $uv$ and treated as a set $\{u, v\}$ (an ordered pair $(u,v)$). We define $\delta(v) := \{e \in E(G) \mid e = uv \textrm{ for some } u \in V\}$.

For matroid-theoretic concepts we follow the notation of Oxley~\cite{oxley}. Familiarity with the definitions and results in that work is assumed.

\subsection{The partial-field axioms}
The following definitions are taken from Semple and Whittle~\cite{SW96}.

\begin{definition}
Let $P$ be a set with distinguished elements called $0$, $1$. Suppose $\cdot$ is a binary operation and $+$ a partial binary operation on $P$. A \emph{partial field} is a quintuple
\begin{align}
  \parf := (P,+, \cdot,0,1)
\end{align}
satisfying the following axioms:
\let\saveenumi\theenumi
\renewcommand{\theenumi}{P\arabic{enumi}}
\begin{enumerate}
    \item \label{ax:mul}$(P\setminus \{0\},\cdot,1)$ is an abelian group.
    \item \label{ax:zero}For all $p \in P$, $p + 0 = p$.
    \item \label{ax:additiveInverse}For all $p \in P$, there is a unique element $q \in P$ such that $p + q = 0$. We denote this element by $-p$.
    \item \label{ax:additiveComm}For all $p, q \in P$, if $p+q$ is defined, then $q+p$ is defined and $p+q=q+p$.
    \item \label{ax:distributivity}For all $p,q,r \in P$, $p\cdot(q+r)$ is defined if and only if $p\cdot q+p\cdot r$ is defined. Then $p\cdot (q+r) = p\cdot q+p\cdot r$.
    \newcounter{savecounter}
    \setcounter{savecounter}{\value{enumi}}
    \item \label{ax:associativelaw}The \emph{associative law} holds for $+$.
\end{enumerate}
\let\theenumi\saveenumi
\end{definition}
If $p,q \in P$ then we abbreviate $p \cdot q$ to $pq$. We write $p+q \doteq r$ if we mean ``the sum of $p$ and $q$ is defined and is equal to $r$''. The group in Axiom~\eqref{ax:mul} is denoted by $\parf^*$, and we write $p \in \parf$ if $p$ is an element of the set $P$ underlying the partial field.

Given a multiset $S = \{p_1, \ldots, p_n\}$ of elements of $P$, a \emph{pre-association} is a vertex-labelled binary tree $T$ with root $r$ such that the leaves are labelled with the elements of $S$ (and each element labels a unique leaf). Moreover, let $v$ be a non-leaf node of $T - r$ with children labelled $u$, $w$. Then $u+w$ must be defined and $v$ is labelled by $u + w$. If $u$, $w$ are the labels of the children of $r$ and $u+w$ is defined, then the labelled tree obtained from $T$ by labeling $r$ with $u+w$ is called an \emph{association} of $S$.

Let $T$ be an association for $S$ with root node $r$, and let $T'$ be a pre-association for the same set (but possibly with completely different tree and labeling). Let $u'$, $w'$ be the labels of the children of the root node of $T'$. Then $T'$ is \emph{compatible} with $T$ if $u'+w' \doteq r$. The \emph{associative law} is the following:
\let\saveenumi\labelenumi
\renewcommand{\labelenumi}{(P\arabic{enumi})}
\begin{enumerate}
  \addtocounter{enumi}{\value{savecounter}}
  \item For every multiset $S$ of elements of $P$ for which some association $T$ exists, every pre-association of $S$ is compatible with $T$.
\end{enumerate}
\let\labelenumi\saveenumi

We say that the expression $p_1 + \cdots + p_n$ is \emph{defined} if there exists a finite multiset $Z$ of the form $\{z_1, -z_1, z_2, -z_2, \ldots, z_k, -z_k\}$ such that there exists an association for $S := \{p_1, \ldots, p_n\} \cup Z$. The value of $p_1+ \cdots + p_n$ is then defined as the value of $r$ for any association $T$ of $S$. Note that this definition differs from the one given by Semple and Whittle. A justification for this modification is given in Appendix~\ref{sec:counterex}.

Partial fields share several basic properties with fields. We use the following implicitly in this  paper:
\begin{proposition}\label{lem:basics}
Let $\parf$ be a partial field. The following statements hold for all $p, q \in \parf$:
\begin{enumerate}
  \item $0 p = 0$;
  \item $p q = 0$ if and only if $p = 0$ or $q = 0$;
  \item $(-1)^2 = 1$;
  \item if $p + q \doteq r$, then $r - q \doteq p$.
\end{enumerate}
\end{proposition}
The proofs are elementary.

\subsection{Partial-field matrices}
Recall that formally, for ordered sets $X$ and $Y$, an $X\times Y$ matrix $A$ with entries in a partial field $\parf$ is a function $A:X\times Y \rightarrow \parf$.
Let $A$ be an $n\times n$ matrix with entries in $\parf$. Then the \emph{determinant} of $A$ is, as always,
\begin{align}
  \det(A) := \sum_{\sigma \in S_n} \sign(\sigma) a_{1\sigma(1)}a_{2\sigma(2)}\cdots a_{n\sigma(n)}.
\end{align}
We say that $\det(A)$ is \emph{defined} if this sum is defined.
\begin{proposition}[{\cite[Proposition 3.1]{SW96}}]\label{prop:matops}Let $\parf$ be a partial field and let $A$ be an $n\times n$ matrix with entries in $\parf$ such that $\det(A)$ is defined.
\begin{enumerate}
\item\label{eq:matopsone} If $B$ is obtained from $A$ by transposition, then $\det(B) \doteq \det(A)$.
\item\label{eq:matopstwo} If $B$ is obtained from $A$ by interchanging a pair of rows, then $\det(B) \doteq -\det(A)$.
\item\label{eq:matopsthree} If $B$ is obtained from $A$ by multiplying a row by a non-zero element $p \in \parf$, then $\det(B) \doteq p \det(A)$.
\item\label{eq:matopsfour} If $B$ is obtained from $A$ by adding two rows whose sum is defined, then $\det(B) \doteq \det(A)$.
\end{enumerate}
\end{proposition}
An $X\times Y$ matrix $A$ with entries in $\parf$ is a \emph{$\parf$-matrix} if $\det(A')$ is defined for every square submatrix $A'$ of $A$. For such a matrix we define the \emph{rank}
\begin{align}
  \rank(A) := \max \{r \mid A \textrm{ has an } r\times r \textrm{ submatrix } A' \textrm{ with } \det(A') \neq 0\}.
\end{align}

Let $A$ be an $X\times Y$ $\parf$-matrix such that $X\cap Y = \emptyset$, and let $x \in X, y \in Y$ be such that $A_{xy} \neq 0$. Then we define $A^{xy}$ to be the $(X\setminus x \cup y) \times (Y\setminus y \cup x)$ matrix with entries
\begin{align}
  (A^{xy})_{uv} = \left\{ \begin{array}{ll}
    A_{xy}^{-1} \quad & \textrm{if } uv = yx\\
    A_{xy}^{-1} A_{xv} & \textrm{if } u = y, v\neq x\\
    -A_{xy}^{-1} A_{uy} & \textrm{if } v = x, u \neq y\\
    A_{uv} - A_{xy}^{-1} A_{uy} A_{xv} & \textrm{otherwise.}
  \end{array}\right.
\end{align}
We say that $A^{xy}$ is obtained from $A$ by \emph{pivoting} over $xy$. In other words, if $X = X' \cup x$, $Y = Y'\cup y$, and
\begin{align}
  A = \kbordermatrix{
          & y & & Y'\\
        x & a & \vrule & b\\ \cline{2-4}
        X' & c & \vrule & D
      },
\end{align}
where $a \in \parf^*$ (i.e. $a \neq 0$), $b$ is a row vector, $c$ a column vector, and $D$ an $X'\times Y'$ matrix, then
\begin{align}
  A^{xy} = \kbordermatrix{
          & x & & Y'\\
        y & a^{-1} & \vrule & a^{-1} b\\ \cline{2-4}
        X' & -a^{-1}c  & \vrule & D - a^{-1}c b
      }.
\end{align}
We refer readers who are unfamiliar with the pivot operation to Oxley~\cite[Page 84; Page 209]{oxley}.

\begin{definition}\label{def:matops}Let $A$ be an $X\times Y$ $\parf$-matrix, such that $X\cap Y = \emptyset$. We say that $A'$ is a \emph{minor} of $A$ (notation: $A' \minorof A$) if $A'$ can be obtained from $A$ by a sequence of the following operations:
\begin{enumerate}
  \item Multiplying the entries of a row or column by an element of $\parf^*$;
  \item Deleting rows or columns;
  \item Permuting rows or columns (and permuting labels accordingly);
  \item Pivoting over a nonzero entry.
\end{enumerate}
\end{definition}
Be aware that in linear algebra a minor of a matrix has a different definition. We use Definition~\ref{def:matops} because of its relation with matroid minors, which will be explained in the next section. For a determinant of a square submatrix we use the word \emph{subdeterminant}.
\begin{proposition}[{\cite[Proposition 3.3]{SW96}}]\label{prop:pmatops}Let $A$ be a $\parf$-matrix. Then $A^T$ is also a $\parf$-matrix. If $A'\minorof A$ then $A'$ is a $\parf$-matrix.
\end{proposition}
If $X'\subseteq X$ and $Y'\subseteq Y$, then we denote by $A[X',Y']$ the submatrix of $A$ obtained by deleting all rows and columns in $X\setminus X'$, $Y\setminus Y'$. If $Z$ is a subset of $X\cup Y$ then we define $A[Z] := A[X\cap Z, Y \cap Z]$. Also, $A-Z := A[X\setminus Z, Y\setminus Z]$. The following observation is used throughout this paper:
\begin{lemma}\label{lem:detpivot}
  Let $A$ be an $X\times Y$ matrix with entries in $\parf$ such that $X\cap Y = \emptyset$ and $|X|=|Y|$. If $\det(A^{xy}-\{x,y\})$ is defined then $\det(A)$ is defined, and
  \begin{align}
  \det(A)= (-1)^s A_{xy} \det(A^{xy}-\{x,y\})
  \end{align}
  for some $s \in \{0,1\}$.
\end{lemma}

Let $A$ be an $X\times Y$ $\parf$-matrix, and let $A'$ be an $X'\times Y'$ $\parf$-matrix. Then $A$ and $A'$ are \emph{isomorphic} if there exist bijections $f:X\rightarrow X'$, $g:Y\rightarrow Y'$ such that for all $x\in X$, $y \in Y$, $A_{xy} = A'_{f(x)g(y)}$.

Let $A$, $A'$ be $X\times Y$ $\parf$-matrices. If $A'$ can be obtained from $A$ by scaling rows and columns by elements from $\parf^*$, then we say that $A$ and $A'$ are \emph{scaling-equivalent}, which we denote by $A\sim A'$.

Let $A$ be an $X\times Y$ $\parf$-matrix such that $X\cap Y = \emptyset$, and let $A'$ be an $X'\times Y'$ $\parf$-matrix such that $X\cup Y = X'\cup Y'$. If $A' \minorof A$ and $A \minorof A'$, then we say that $A$ and $A'$ are \emph{strongly equivalent}, which we denote by $A' \approx A$. If $\phi(A') \approx A$ for some partial field automorphism $\phi$ (see below for a definition), then we say $A'$ and $A$ are \emph{equivalent}.
\subsection{Partial-field matroids}
Let $A$ be an $r\times E$ $\parf$-matrix of rank $r$. We define the set
\begin{align}
  \mathcal{B}_A := \{B \subseteq E \mid |B| = r, \det(A[r,B])\neq 0\}.
\end{align}
\begin{theorem}[{\cite[Theorem 3.6]{SW96}}]$\mathcal{B}_A$ is the set of bases of a matroid.
\end{theorem}
We denote this matroid by $M[A] = (E,\mathcal{B}_A)$.
Conversely, let $M$ be a matroid. If there exists a $\parf$-matrix $A$ such that $M = M[A]$, then we say that $M$ is $\parf$-representable. These matroids share many properties of representable matroids.
\begin{lemma}[{\cite[Proposition 4.1]{SW96}}]
  Let $A$ be an $r\times E$ $\parf$-matrix, and $B$ a basis of $M[A]$. Then there exists a $\parf$-matrix $A'$ such that $M[A'] = M[A]$ and $A'[r,B]$ is an identity matrix.
\end{lemma}
Conversely, let $A$ be an $X\times Y$ matrix with entries in $\parf$, such that $X\cap Y = \emptyset$. Let $A'$ be the $X\times (X\cup Y)$ matrix $A' = [I|A]$, where $I$ is an $X\times X$ identity matrix. For all $X' \subseteq X\cup Y$ with $|X'| = |X|$ we have $\det(A'[X,X']) = \pm \det(A[X\setminus X', Y\cap X'])$. Hence $A'$ is a $\parf$-matrix if and only if $A$ is a $\parf$-matrix. We say that $M = M[I|A]$ is the matroid \emph{associated with} $A$, and that $[I|A]$ is an \emph{$X$-representation} of $M$ for basis $X$.

If $N$ is a minor of a matroid $M$, say $N = M\delete S \contract T$, then a $B$-representation \emph{displays} $N$ if $B\cap T = T$ and $B\cap S = \emptyset$; then $N = M[I'|A']$, where $A' = A-S-T$. Likewise we say that $A$ \emph{displays} $A'$ if $A' = A-U$ for some $U \subseteq X \cup Y$.
\begin{lemma}\label{lem:matrixmatroid}
  If $M = M[I|A]$, then $N\minorof M$ if and only if $N\cong M[I'|A']$ for some $A'\minorof A$.
\end{lemma}

\subsection{Partial-field homomorphisms}
A function $\phi: \parf_1 \rightarrow \parf_2$ is a \emph{homomorphism} if, for all $p, q \in \parf_1$, $\phi(pq) = \phi(p)\phi(q)$ and, when $p+q$ is defined, then $\phi(p)+\phi(q) \doteq \phi(p+q)$. A homomorphism is \emph{trivial} if its kernel is equal to $\parf_1$. This happens if and only if $\phi(1) = 0$.
\begin{proposition}[{\cite[Proposition 5.1]{SW96}}]\label{prop:hom}Let $\parf_1$, $\parf_2$ be partial fields and let $\phi: \parf_1 \rightarrow \parf_2$ be a homomorphism. Let $A$ be a $\parf_1$-matrix. Then
\begin{enumerate}
\item $\phi(A)$ is a $\parf_2$-matrix.
\item If $A$ is square and $\det(A) = 0$ then $\det(\phi(A)) = 0$.
\item If $A$ is square and $\phi$ is nontrivial then $\det(A) = 0$ if and only if $\det(\phi(A)) = 0$.
\end{enumerate}
\end{proposition}
This leads to the following easy corollary:
\begin{corollary}[{\cite[Corollary 5.3]{SW96}}]\label{cor:hom}Let $\parf_1$ and $\parf_2$ be partial fields and let $\phi: \parf_1 \rightarrow \parf_2$ be a nontrivial homomorphism. If $A$ is a $\parf_1$-matrix then $M[\phi(A)] = M[A]$. It follows that, if $M$ is a $\parf_1$-representable matroid, then $M$ is also $\parf_2$-representable.
\end{corollary}
A partial field isomorphism $\phi:\parf_1\rightarrow\parf_2$ is a bijective homomorphism with the additional property that $\phi(p+q)$ is defined if and only if $p+q$ is defined. If $\parf_1$ and $\parf_2$ are isomorphic then we denote this by $\parf_1 \cong \parf_2$. A partial field automorphism is an isomorphism $\phi:\parf\rightarrow\parf$.

\subsection{Constructions}\label{ssec:sources}
For a general partial field the associative law is hard to wield. Semple and Whittle get around this difficulty by constructing partial fields as restrictions of bigger partial fields, starting their construction with a field. Recall that $\parf^*$ is the multiplicative group of $\parf$, and for $S\subseteq \parf^*$, $\langle S \rangle$ is the subgroup generated by $S$.
\begin{definition}\label{def:subpf}
  Let $\parf$ be a partial field, and let $S$ be a set of elements of $\parf^*$. Then
\begin{align}
\parf[S]:= (\langle S\cup -1\rangle\cup 0, 0,1,+,\cdot),
\end{align}
where multiplication and addition are the restriction of the operations in $\parf$, i.e. $p+q$ is defined only if $p+q\doteq r$ in $\parf$ and $r \in \langle S\cup-1 \rangle\cup 0$.
\end{definition}
\begin{proposition}[{\cite[Proposition 2.2]{SW96}}]\label{prop:subpf}$\parf[S]$ is a partial field.
\end{proposition}
We need $-1 \in \parf[S]$ to ensure that $1$ has an additive inverse.

Instead of constructing a partial field as the restriction of a field, one can also take a ring as starting structure.
\begin{definition}\label{prop:pffromring}
Let $\ring$ be a commutative ring, and let $S$ be a subset of $\ring^*$. Then
\begin{align}
  \parf(\ring,S) := (\langle S\cup -1\rangle\cup 0, 0, 1, + ,\cdot),
\end{align}
where multiplication and addition are the restriction of the operations in $\ring$, i.e. $p+q$ is defined only if the resulting element of $\ring$ is again in $\langle S \cup -1\rangle\cup 0$.
\end{definition}
\begin{proposition}$\parf(\ring, S)$ is a partial field.
\end{proposition}
\begin{proof}
First remark that $1 \in \parf$ and that $-1$ is invertible in $\ring$. The other axioms are then inherited from the corresponding ring axioms.
\end{proof}
In fact, Proposition~\ref{prop:subpf} is a special case of this result. To see this we need to find a suitable ring. The the following theorem provides such a ring:
\begin{theorem}[{Vertigan}]\label{thm:pfinring}If $\parf$ is a partial field, then there exist a ring $\ring$ and a set $S \subseteq\ring^*$ such that $\parf \cong \parf(\ring,S)$.
\end{theorem}
We present a proof of this theorem in Section~\ref{sec:rings}.
A third source of partial fields is the following. If $\parf_1$, $\parf_2$ are partial fields, then we define the \emph{direct product}
\begin{align}
  \parf_1\otimes\parf_2 := (P, +,\cdot,(0,0),(1,1)),
\end{align}
where
\begin{align}
  P = \{(p_1,p_2) \in \parf_1\times \parf_2 \mid p_1 \neq 0 \textrm{ if and only if } p_2 \neq 0\}
\end{align}
and addition and multiplication are defined componentwise, i.e. $(p_1,p_2) + (q_1,q_2)\doteq (p_1+q_1, p_2+q_2)$ if and only if both $p_1+q_1$ and $p_2+q_2$ are defined and $p_1 + q_1 = 0$ if and only if $p_2 + q_2 = 0$.
\begin{lemma}
  $\parf_1\otimes \parf_2$ is a partial field.
\end{lemma}
\begin{proof}
  This follows from an application of Proposition~\ref{prop:pffromring}:
if $\parf_i = \parf(\ring_i,S_i)$ then $\parf_1 \otimes \parf_2 = \parf(\ring_1\times \ring_2, S_1\times S_2)$.
\end{proof}
Suppose $\parf$, $\parf_1$, $\parf_2$ are partial fields such that there exist homomorphisms $\phi_1:\parf\rightarrow \parf_1$ and $\phi_2:\parf\rightarrow\parf_2$. Then we define $\phi_1\otimes \phi_2: \parf \rightarrow \parf_1\otimes\parf_2$ by $(\phi_1\otimes\phi_2)(p) := (\phi_1(p),\phi_2(p))$.
\begin{lemma}$\phi_1\otimes\phi_2$ is a partial field homomorphism.
\end{lemma}
The proof is straightforward and therefore omitted.

Let $X$, $Y$ be finite, disjoint sets, let $A_1$ be an $X\times Y$ $\parf_1$-matrix, and let $A_2$ be an $X\times Y$ $\parf_2$-matrix. Let $A := A_1\otimes A_2$ be the $X\times Y$ matrix such that $A_{uv} = ((A_1)_{uv},(A_2)_{uv})$.
\begin{lemma}If $A_1$ is a $\parf_1$-matrix, $A_2$ is a $\parf_2$-matrix, and $M[I|A_1] = M[I|A_2]$ then $A_1\otimes A_2$ is a $\parf_1\otimes\parf_2$-matrix and $M[I|A_1\otimes A_2] = M[I|A_1]$.
\end{lemma}
\begin{proof}
  Let $X'\subseteq X$, $Y' \subseteq Y$ such that $A' := A[X',Y']$ is a square submatrix of $A$. Since $M[I|A_1] = M[I|A_2]$, $\det(A_1[X',Y']) = 0$ if and only if $\det(A_2[X',Y']) = 0$. This holds for all $1\times 1$ submatrices as well, so all entries of $A$ are from $\parf_1\otimes\parf_2$. By Lemma~\ref{lem:detpivot}, a determinant can be computed by a sequence of pivots. It follows that $\det(A')$ is defined, which completes the proof.
\end{proof}
The following corollary plays a central role in this paper.
\begin{corollary}\label{thm:listoffields}Let $M$ be a matroid. $M$ is representable over each of $\parf_1, \ldots, \parf_k$ if and only if it is representable over the partial field
\begin{align}
  \parf := \parf_1 \otimes \cdots \otimes \parf_k.
\end{align}
\end{corollary}

\subsection{Cross ratios and fundamental elements}\label{ssec:funds}
Let $B = \left[\begin{smallmatrix}p & q\\r & s\end{smallmatrix}\right]$ be a $\parf$-matrix with $ps \neq 0$. We define the \emph{cross ratio} of $B$ as
\begin{align}
  \cra(B) := \frac{qr}{ps}.
\end{align}
The motivation for this name comes from projective geometry. If $\cra(B) \not \in \{0,1\}$ then the matroid $M[I|B]$ is the four-point line. In projective geometry the cross ratio is a number defined for any ordered set of four collinear points. It is invariant under projective transformations. For a fixed set of points this number can take six different values, depending on the order.

Let $A$ be an $X\times Y$ $\parf$-matrix. We define the \emph{cross ratios} of $A$ as the set
\begin{align}
    \crat(A) := \left\{\cra\left(\left[\begin{smallmatrix}1 & 1\\p & 1\end{smallmatrix}\right]\right) \mid \left[\begin{smallmatrix}1 & 1\\p & 1\end{smallmatrix}\right] \minorof A \right\}.
\end{align}
The following is obvious from the definition:
\begin{lemma}\label{lem:cratminor}
  If $A' \minorof A$ then $\crat(A') \subseteq \crat(A)$.
\end{lemma}

Note that $\det\left(\left[\begin{smallmatrix}1 & 1\\p & 1\end{smallmatrix}\right]\right) = 1-p$. This prompts the following definition. An element $p \in \parf$ is called \emph{fundamental} if $1-p \in \parf$. As remarked by Semple~\cite{Sem97}, $p+q$ is defined if and only if $p^{-1} (p+q) = 1 - (-q/p)$ is defined. For most partial fields that we consider, the equation $1-p = q$ has only finitely many solutions. This is convenient if one wants to compute in partial fields (cf. Hlin\v en\'y~\cite{Hli04}). We denote the set of fundamental elements of $\parf$ by $\fun(\parf)$.

Suppose $F\subseteq \fun(\parf)$. We define the \emph{associates} of $F$ as
\begin{align}
  \assoc F := \bigcup_{p\in F}\crat\left( \left[\begin{smallmatrix}
    1 & 1\\
    p & 1
  \end{smallmatrix}\right]
  \right).
\end{align}
We have
\begin{proposition}
  $\assoc\{p\}\subseteq \fun(\parf)$.
\end{proposition}
The following lemma gives a complete description of the structure of $\assoc\{p\}$.
\begin{lemma}If $p \in \{0,1\}$ then $\assoc\{p\} = \{0,1\}$. If $p \in \fun(\parf)\setminus\{0,1\}$ then
\begin{align}
  \assoc\{p\} = \big\{p, 1-p,\frac{1}{1-p},\frac{p}{p-1},\frac{p-1}{p},\frac{1}{p}\big\}.\label{eq:associates}
\end{align}
\end{lemma}
The proof consists of a straightforward enumeration.
By Lemma~\ref{lem:cratminor}, $\assoc\{p\} \subseteq \crat(A)$ for every $p \in \crat(A)$.

\subsection{Normalization}
Let $M$ be a rank-$r$ matroid with ground set $E$, and let $B$ be a basis of $M$. Let $G=G(M,B)$ be the bipartite graph with vertices $V(G) = B \cup (E\setminus B)$ and edges $E(G) = \{xy \in B\times (E\setminus B) \mid (B\setminus x) \cup y\in \mathcal{B}\}$. For each $y \in E\setminus B$ there is a unique matroid circuit $C_{B,y} \subseteq B\cup y$, the $B$-\emph{fundamental circuit} of $y$.
\begin{lemma}\label{lem:bipconn}
  Let $M$ be a matroid, and $B$ a basis of $M$.
  \begin{enumerate}
    \item $xy \in E(G)$ if and only if $x \in C_{B,y}$.
    \item\label{bip:conn} $M$ is connected if and only if $G(M,B)$ is connected.
    \item\label{bip:threeconn} If $M$ is 3-connected, then $G(M,B)$ is 2-connected.
  \end{enumerate}
\end{lemma}
\begin{proof}
  This follows from consideration of the $B$-fundamental-circuit incidence matrix. See, for example, Oxley~\cite[Section 6.4]{oxley}.
\end{proof}

Let $A$ be an $X\times Y$ matrix, such that $X\cap Y = \emptyset$. With $A$ we associate a bipartite graph $G=\bip(A)$ with vertices $V(G) = X \cup Y$ and edges $E(G) = \{xy \in X\times Y \mid A_{xy} \neq 0\}$. Recall that $\sim$ denotes scaling-equivalence.
\begin{lemma}\label{lem:bipproperties}
  Let $\parf$ be a partial field. Suppose $M=M[I|A]$.
  \begin{enumerate}
    \item \label{eq:bipnonzeroentries}$G(M,X) = \bip(A)$.
    \item \label{eq:bipscaling}Let $T$ be a spanning forest of $\bip(A)$ with edges $e_1, \ldots, e_k$. Let $p_1, \ldots, p_k \in \parf^*$. Then there exists a matrix $A'\sim A$ such that $A'_{e_i} = p_i$.
  \end{enumerate}
\end{lemma}
The proof of the corresponding theorem in Oxley~\cite[Theorem 6.4.7]{oxley} generalizes directly to partial fields.

Let $A$ be a matrix and $T$ a spanning forest for $\bip(A)$. We say that $A$ is \emph{$T$-normalized} if $A_{xy} = 1$ for all $xy \in T$. By the lemma there is always an $A'\sim A$ that is $T$-normalized. We say that $A$ is \emph{normalized} if it is $T$-normalized for some spanning forest $T$, the \emph{normalizing} spanning forest.

The following definitions are needed for the statement and proof of Theorem~\ref{thm:lift}. As usual, a \emph{walk} in a graph $G = (V,E)$ is a sequence $W = (v_0, \ldots, v_n)$ of vertices such that $v_iv_{i+1} \in E$ for all $i \in \{0, \ldots, n-1\}$. If $v_n = v_0$ and $v_i \neq v_j$ for all $0 \leq i < j < n$ then we say that $W$ is a \emph{cycle}.
\begin{definition}
  Let $A$ be an $X\times Y$ matrix with entries in a partial field $\parf$, such that $X\cap Y = \emptyset$. The \emph{signature} of $A$ is the function $\sigma_A:(X\times Y) \cup (Y\times X) \rightarrow \parf$ defined by
  \begin{align}
    \sigma_A(vw) := \left\{ \begin{array}{cl} A_{vw}\quad & \textrm{if } v \in X, w \in Y\\
                                             1/A_{vw}\quad & \textrm{if } v \in Y, w \in X.
                           \end{array}\right.
  \end{align}
  If $C = (v_0, v_1, \ldots, v_{2n-1},v_{2n})$ is a cycle of $\bip(A)$ then we define
  \begin{align}
    \sigma_A(C) := (-1)^{|V(C)|/2}\prod_{i=0}^{2n-1} \sigma_A(v_iv_{i+1}).
  \end{align}
\end{definition}
Observe that the signature of a cycle does not depend on the choice of $v_0$. If $C'$ is the cycle $(v_{2n}, v_{2n-1},\ldots, v_1, v_0)$ then $\sigma_A(C') = 1/\sigma_A(C)$. If $A$ a $\parf$-matrix such that $\bip(A)$ is a cycle, then $M[I|A]$ is a wheel if the signature equals $1$, and a whirl otherwise.

The proof of the following lemma is straightforward. The last property exhibits a close connection between the signature and determinants. Recall that $A^{xy}$ is the matrix obtained from $A$ by pivoting over $xy$.
\begin{lemma}\label{lem:signature}
  Let $A$ be an $X\times Y$ matrix with entries from a partial field $\parf$, such that $X\cap Y = \emptyset$.
  \begin{enumerate}
    \item\label{sgn:scale} If $A'\sim A$ then $\sigma_{A'}(C)= \sigma_A(C)$ for all cycles $C$ in $\bip(A)$.
    \item\label{sgn:pivot} Let $C = (v_0, \ldots, v_{2n})$ be an induced cycle of $\bip(A)$ with $v_0 \in X$ and $n \geq 3$. Suppose $A' := A^{v_0v_1}$ is such that all entries are defined. Then $C' = (v_2, v_3, \ldots, v_{2n-1},v_2)$ is an induced cycle of $\bip(A')$ and $\sigma_{A'}(C') = \sigma_A(C)$.
    \item\label{sgn:det} Let $C = (v_0, \ldots, v_{2n})$ be an induced cycle of $\bip(A)$. If $A'$ is obtained from $A$ by scaling rows and columns so that $A'_{v_iv_{i+1}} = 1$ for all $i > 0$, then $A'_{v_0v_1} = (-1)^{|V(C)|/2} \sigma_A(C)$ and $\det(A[V(C)]) = 1-\sigma_A(C)$.
  \end{enumerate}
\end{lemma}
\begin{corollary}\label{cor:signaturefund}
  Let $A$ be an $X\times Y$ $\parf$-matrix. If $C$ is an induced cycle of $\bip(A)$ then $\sigma_A(C) \in \crat(A) \subseteq \fun(\parf)$.
\end{corollary}

\subsection{Examples}
We can now give a very short proof of Theorem~\ref{thm:reg}. First we restate it using our new terminology. We define the \emph{regular} partial field
\begin{align}
  \reg := \parf(\Q,\emptyset).
\end{align}
It has just three elements: $\{-1,0,1\}$. Clearly a $\reg$-matrix is a totally unimodular matrix.
\begin{theorem}[Tutte \cite{Tut65}]
Let $M$ be a matroid. The following are equivalent:
\begin{enumerate}
  \item \label{eq:regtwothree}$M$ is representable over $\GF(2)\otimes\GF(3)$;
  \item \label{eq:regTU} $M$ is $\reg$-representable.
  \item \label{eq:regevery}$M$ is representable over every partial field.
\end{enumerate}
\end{theorem}
\begin{proof}
  Every partial field $\parf$ contains a multiplicative identity and, by Axiom~\eqref{ax:additiveInverse}, an element $-1$. Therefore there exists a nontrivial homomorphism $\phi:\reg \rightarrow \parf$, which proves \eqref{eq:regTU}$\Rightarrow$\eqref{eq:regevery}. The partial field $\GF(2)\otimes\GF(3)$ has fundamental elements $\{(0,0),(1,1)\}$. We have an obvious homomorphism $\phi':\GF(2)\otimes\GF(3) \rightarrow \reg$, which proves \eqref{eq:regtwothree}$\Rightarrow$\eqref{eq:regTU}. Finally, \eqref{eq:regevery}$\Rightarrow$\eqref{eq:regtwothree} is trivial.
\end{proof}
We define the \emph{sixth roots of unity} partial field $\psru := \parf(\C,\zeta)$, where $\zeta$ is a root of $x^2-x+1=0$, i.e. $\zeta$ is a primitive sixth root of unity. Whittle proved the following theorem:
\begin{theorem}[Whittle~\cite{Whi97}]\label{thm:sru}
Let $M$ be a matroid. The following are equivalent:
\begin{enumerate}
  \item \label{eq:srutwothree}$M$ is representable over $\GF(3)\otimes\GF(4)$;
  \item \label{eq:sruTU} $M$ is $\psru$-representable;
  \item \label{eq:sruevery} $M$ is representable over $\GF(3)$, over $\GF(p^2)$ for all primes $p$, and over $\GF(p)$ when $p \equiv 1 \mod 3$.
\end{enumerate}
\end{theorem}
\begin{proof}
  Note that $\psru$ is finite, with $\fun(\psru) = \{0,1,\zeta, 1-\zeta\}$. Let $\phi:\psru\rightarrow\GF(3)\otimes\GF(4)$ be determined by $\phi(\zeta) = (-1,\omega)$, where $\omega \in \GF(4)\setminus \{0,1\}$ is a generator of $\GF(4)^*$. Then $\phi$ is a bijective homomorphism, which proves \eqref{eq:srutwothree}$\Leftrightarrow$\eqref{eq:sruTU}.

  That \eqref{eq:sruevery} implies \eqref{eq:srutwothree} is again trivial. We will use results from algebraic number theory to prove \eqref{eq:sruTU}$\Rightarrow$\eqref{eq:sruevery}. See, for example, Stewart and Tall~\cite{ST87} for the necessary background.
  For \eqref{eq:sruTU}$\Rightarrow$\eqref{eq:sruevery}, remark that $\psru^*$ is the group of units of $\Z[\zeta]$, the ring of integers of the algebraic number field $\Q(\zeta) = \Q(\sqrt{-3})$. If $I$ is a maximal ideal then $\Z[\zeta]/I$ is a finite field. We find the values $q = p^m$ for which there exists a prime ideal $I$ with norm $N(I) := |\Z[\zeta]/I| = q$. If $I$ is a principal ideal, i.e. $I = (a+b\sqrt{-3})\Z[\zeta]$ with $a,b \in \frac{1}{2}\Z$, then $N(I) = a^2 + 3 b^2$.

  Suppose $I = (\sqrt{-3})\Z[\zeta]$. Then $N(I) = 3$ which is prime, so $\Z[\zeta]/I \cong \GF(3)$. This gives a ring homomorphism $\phi:\Z[\zeta]\rightarrow\GF(3)$. Suppose $I = p\Z[\zeta]$. Then $N(p\Z[\zeta]) = p^2$. Either $I$ is prime, in which case $\Z[\zeta]/I \cong \GF(p^2)$, or $I$ splits and there exists a prime ideal $J$ with $\Z[\zeta]/J \cong \GF(p)$. A well-known result in number theory (see e.g. Hardy and Wright~\cite[{Theorem 255}]{HW54}) states that $I$ splits if and only if $p \equiv 1 \mod 3$.
\end{proof}

Whittle gave characterizations for several other classes of matroids. However, the proofs of these are more complicated, because the partial fields involved are no longer isomorphic. In the next section we develop a general tool to overcome this difficulty.

\section{The lift theorem}\label{sec:lifts}
Let $\parf$, $\hat\parf$ be partial fields and let $\phi:\hat\parf\rightarrow\parf$ be a homomorphism. Let $A$ be an $X\times Y$ $\parf$-matrix. In what follows we would like to construct an $X\times Y$ $\hat\parf$-matrix $\hat A$ such that $\phi(\hat A) = A$, even in the absence of a partial field homomorphism $\parf\rightarrow\hat\parf$. To that end we make the following definitions. Recall that $\fun(\parf)$ is the set of fundamental elements of a partial field.
\begin{definition}
  Let $\parf$, $\hat \parf$ be partial fields, and let $\phi:\hat\parf\rightarrow\parf$ be a partial field homomorphism. A \emph{lifting function} for $\phi$ is a function $^\uparrow:\fun(\parf)\rightarrow\hat\parf$ such that for all $p,q\in\fun(\parf)$:
\begin{itemize}
  \item $\phi(p^\uparrow) = p$;
  \item if $p+q \doteq 1$ then $p^\uparrow + q^\uparrow \doteq 1$;
  \item if $p\cdot q = 1$ then $p^\uparrow\cdot q^\uparrow = 1$.
\end{itemize}
\end{definition}
Hence a lifting function maps $\assoc \{p\}$ to $\assoc\{p^\uparrow\}$ for all $p \in \fun(\parf)$.
\begin{definition}\label{def:locallift}
  Let $\parf$, $\hat\parf$ be two partial fields, let $\phi: \hat\parf \rightarrow \parf$ be a homomorphism, and let $^\uparrow:\fun(\parf)\rightarrow\hat\parf$ be a lifting function for $\phi$. Let $A$ be an $X\times Y$ $\parf$-matrix. An $X\times Y$ matrix $\hat A$ is a \emph{local $^\uparrow$-lift} of $A$ if
  \begin{enumerate}
    \item\label{loclift:one} $\phi(\hat A) \sim A$;
    \item\label{loclift:two} $\hat A$ is an $X\times Y$ $\hat\parf$-matrix;
    \item\label{loclift:three} for every induced cycle $C$ of $\bip(A)$ we have
      \begin{align}\label{eq:cyclegood}
        \sigma_A(C)^\uparrow = \sigma_{\hat A}(C).
      \end{align}
  \end{enumerate}
\end{definition}
First we show that, if a local $^\uparrow$-lift exists, it is unique up to scaling.
\begin{lemma}\label{lem:uniqA2}
  Let $\parf$, $\hat\parf$ be two partial fields, let $\phi: \hat\parf \rightarrow \parf$ be a homomorphism, and let $^\uparrow:\fun(\parf)\rightarrow\hat\parf$ be a lifting function for $\phi$. Let $A$ be an $X\times Y$ $\parf$-matrix, and suppose $\hat A_1$, $\hat A_2$ are local $^\uparrow$-lifts of $A$. Then $\hat A_1\sim \hat A_2$.
\end{lemma}
\begin{proof}
  Suppose the lemma is false and let $A$, $\hat A_1$, $\hat A_2$ form a counterexample. Let $T$ be a spanning forest of $\bip(A)$ and rescale $\hat A_1$, $\hat A_2$ so that they are $T$-normalized. Let $H$ be the subgraph of $\bip(A)$ consisting of all edges $x'y'$ such that $(\hat A_1)_{x'y'} = (\hat A_2)_{x'y'}$. Let $xy$ be an edge not in $H$ such that the minimum length of an $x-y$ path $P$ in $H$ is minimal. Then $C:=P\cup xy$ is an induced cycle of $\bip(A)$. We have
  \begin{align}
    \sigma_A(C)^\uparrow = \sigma_{\hat A_1}(C) = \sigma_{\hat A_2}(C).
  \end{align}
  But this is only possible if $(\hat A_1)_{xy}=(\hat A_2)_{xy}$, a contradiction.
\end{proof}
It is straightforward to turn this proof into an algorithm that constructs a matrix $\hat A$ satisfying \eqref{loclift:one} and \eqref{loclift:three} for a subset of the cycles such that, if $A$ has a local $^\uparrow$-lift, $\hat A$ is one.

If $\hat A$ is a local lift of $A$, and $A_{xy} \neq 0$, then $\phi(\hat A^{xy}) = A^{xy}$. However, $\hat A^{xy}$ may not be a local lift of $A^{xy}$, since \ref{def:locallift}\eqref{loclift:three} may not hold. This could occur if $\hat \parf$ has more fundamental elements than $\parf$. Next we define a stronger notion of lift, which commutes with pivoting. Recall that we write $A\approx A'$ if $A'$ can be obtained from $A$ by pivoting and scaling.
\begin{definition}
  Let $\parf$, $\hat\parf$ be two partial fields, let $\phi: \hat\parf \rightarrow \parf$ be a homomorphism, and let $^\uparrow:\fun(\parf)\rightarrow\hat\parf$ be a lifting function for $\phi$. A matrix $\hat A$ is a \emph{global $^\uparrow$-lift} of $\phi(\hat A)$ if $\hat A'$ is a local $^\uparrow$-lift of $\phi(\hat A')$ for all $\hat A' \approx \hat A$.
\end{definition}
We now have all ingredients to state the main theorem.
\begin{theorem}[Lift Theorem]\label{thm:lift}Let $\parf$, $\hat\parf$ be two partial fields, let $\phi: \hat\parf \rightarrow \parf$ be a homomorphism, and let $^\uparrow:\fun(\parf)\rightarrow\hat\parf$ be a lifting function for $\phi$.
Let $A$ be an $X\times Y$ $\parf$-matrix. Then exactly one of the following is true:
\begin{enumerate}
  \item \label{enum:propthree}$A$ has a global $^\uparrow$-lift.
  \item \label{enum:propfour} $A$ has a minor $B$ such that
  \begin{enumerate}
    \item $B$ has no local $^\uparrow$-lift;
    \item $B$ or $B^T$ equals
    \begin{align}\label{eq:mincounterexmats}
       \begin{bmatrix}
         0 & 1 & 1 & 1\\
         1 & 0 & 1 & 1\\
         1 & 1 & 0 & 1
       \end{bmatrix} \textrm{ or }
       \begin{bmatrix} 1 & 1 & 1\\
                       1 & p & q
       \end{bmatrix}
    \end{align}
    for some distinct $p, q \in \fun(\parf)\setminus\{0,1\}$.
  \end{enumerate}
\end{enumerate}
\end{theorem}

The matroids $M[I|B]$, where $B$ is as in \eqref{eq:mincounterexmats}, are well-known, and often crop up in matroid theory. They are the fano matroid, $F_7$, the non-fano matroid, $F_7^-$, the five-point line, $U_{2,5}$, and their duals. The fano matroid is an excluded minor for all fields that do not have characteristic $2$.

In the proof of the theorem we use techniques similar to those found in, for example, \cite{Ger89,Tru92,LS99}. In fact, Theorem~\ref{thm:lift} generalizes Gerards'~\cite{Ger89} proof of the excluded-minor characterization for regular matroids. First we prove a graph-theoretic lemma.
\begin{lemma}\label{lem:biptree}
  Let $G = (V,E)$ be a 2-connected bipartite graph with bipartition $(U, W)$. Then either $G$ is a cycle or there exists a spanning tree of $G$ with set of leaves $L$, such that $|L|\geq 3$ and $L \cap U \neq \emptyset$, $L \cap W \neq \emptyset$.
\end{lemma}
\begin{proof}
  Suppose $G$ is a counterexample. Since $G$ is not a cycle, $G$ has a vertex $v$ of degree at least $3$. Let $w_1, w_2, w_3$ be neighbours of $v$, and let $v'$ be a neighbour of $w_1$ other than $v$. Then $(\{v, v', w_1, w_2, w_3\}, \{vw_1, vw_2, vw_3,$ $v'w_1\})$ has $3$ leaves, not all in the same vertex class.

  Now let $T' \subset G$ be a tree with at least three leaves, not all in the same vertex class, such that $V(T')$ is maximal. Let $v \in V(G)\setminus V(T')$. By Menger's Theorem there exist two internally vertex-disjoint $v-T'$ paths $P_1, P_2$. Choose an edge $e \in P_1 \cup P_2$ as follows. If one of the end vertices of $P_1 \cup P_2$ is the unique leaf in $U$ or in $W$, choose $e$ equal to the edge incident with this vertex. Otherwise choose $e$ arbitrarily. Then $(T' \cup P_1 \cup P_2) \setminus e$ is again a tree with the required property. Indeed: adding $P_1$ and $P_2$ to $T'$ destroys at most two leaves. However, deleting $e$ creates equally many leaves again, and if there are two such new leaves, then there is one in each of $U$ and $W$. Note that $T'$ has a third leaf, which remains unaffected by this construction. But this contradicts our initial choice of $T'$, and the proof is complete.
\end{proof}

We also need the following lemma. Semple and Whittle \cite{SW96} proved that the 2-sum of two $\parf$-matrices is again a $\parf$-matrix. We need something slightly stronger.
\begin{lemma}\label{lem:twosum}
  Let $A$ be a $\parf$-matrix, and $(X_1,X_2)$, $(Y_1,Y_2)$ partitions of $X$ and $Y$ such that
  \begin{align}
    A = \kbordermatrix{ & Y_1 & & Y_2\\
                       X_1 & A_1'& \vrule & a_1a_2\\ \cline{2-4}
                       X_2 & 0 & \vrule   & A_2'
         },
  \end{align}
  where $A_1'$, $A_2'$ are submatrices, $a_1$ is a column vector, and $a_2$ is a row vector. If both
  \begin{align}
    A_1 := \begin{bmatrix} A_1' & a_1\\
                             0 & 1
           \end{bmatrix} \textrm{ and } A_2 := \begin{bmatrix}1 & a_2\\ 0 & A_2'\end{bmatrix}
  \end{align}
  have a global $^\uparrow$-lift then $A$ has a global $^\uparrow$-lift.
\end{lemma}
The following proof sketch omits some details, but the remaining difficulties are purely notational.
\begin{proof}[Sketch of proof]
  Let $A$, $A_1$, $A_2$ be as in the lemma, and let $\hat A_1$, $\hat A_2$ be global $^\uparrow$-lifts of $A_1$, $A_2$. We define
  \begin{align}
    \hat A := \kbordermatrix{ & Y_1 & Y_2\\
                       X_1 & \hat A_1' & \hat a_1\hat a_2\\
                       X_2 & 0    & \hat A_2'
         }.
  \end{align}
  By Lemma~\ref{lem:detpivot} every subdeterminant of $\hat A$ is of the form $\pm\det(\hat D_1)\det(\hat D_2)$, where $\hat D_1 \minorof \hat A_1$, and $\hat D_2 \minorof \hat A_2$, from which it follows easily that $\hat A$ is a local lift of $A$. Pick an $x\in X$, $y\in Y$ with $A_{xy} \neq 0$. Then $A^{xy}$ has a minor equivalent to $A_1$ (up to relabelling of rows and columns) and a minor equivalent to $A_2$ (up to relabelling of rows and columns). Moreover $A^{xy}$ can be obtained from these minors in the same way $A$ was obtained from $A_1$ and $A_2$. Therefore $\hat A^{xy}$ must be a local lift of $A^{xy}$. It follows that $A$ has a global lift.
\end{proof}

\begin{proof}[Proof of Theorem \ref{thm:lift}]
\addtocounter{theorem}{-2}
  \eqref{enum:propthree} and \eqref{enum:propfour} cannot hold simultaneously. Suppose the theorem fails for partial fields $\parf$, $\hat\parf$ with homomorphism $\phi$ and lifting function $^\uparrow$. Then there exists a matrix $A$ for which neither \eqref{enum:propthree} nor \eqref{enum:propfour} holds.
  \begin{aclaim}\label{cl:lifttwoconn}
    If $A$ is a counterexample to the theorem with $|X|+|Y|$ minimal then $\bip(A)$ is 2-connected.
  \end{aclaim}
  \begin{subproof}
    If $\bip(A)$ is not connected then one of the components of $A$ has no local $^\uparrow$-lift, contradicting the minimality of $|X|+|Y|$. If $\bip(A)$ has a cut vertex then $A$ is of the form of Lemma~\ref{lem:twosum} with one of $a_1$, $a_2$ having exactly one nonzero entry. Again, the minimality of $|X|+|Y|$ gives a contradiction.
  \end{subproof}
  A pair $(A, \{e,f,g\})$, where $A$ is an $X\times Y$ $\parf$-matrix such that $X\cap Y = \emptyset$, and $\{e,f,g\}\subseteq X\cup Y$, is called a \emph{bad pair} if
  \begin{enumerate}
    \item $A$ is a counterexample to the theorem with $|X|+|Y|$ minimal;
    \item There exists a spanning tree $T$ of $\bip(A)$ such that $\{e,f,g\}$ are leaves of $T$;
    \item $e,f \in X$ and $g \in Y$.
  \end{enumerate}
  \begin{aclaim}\label{cl:minorminimal}
    If $(A, \{e,f,g\})$ is a bad pair then there exists a matrix $\hat A$ such that $\hat A-U$ is a global lift of $A-U$ for all $U$ such that $U\cap \{e,f,g\}\neq\emptyset$.
  \end{aclaim}
  \begin{subproof}
    Without loss of generality $A$ is $T$-normalized for a tree $T$ in which $e,f,g$ are leaves. Note that $T - U$ is a spanning tree of $A-U$ for all nonempty $U\subseteq \{e,f,g\}$. By Lemma~\ref{lem:uniqA2} there exists a unique $T-U$-normalized global $^\uparrow$-lift $\hat{A-U}$ for $A-U$. Again by Lemma~\ref{lem:uniqA2} and our choice of $T$, if $v \in \{e,f,g\}\setminus U$, then $\hat{A-U}-v = \hat{A-U-v}$. It follows that there is a unique matrix $\hat A$ such that $\hat A -U = \hat{A-U}$ for all nonempty $U\subseteq \{e,f,g\}$.
  \end{subproof}
  We say that $\hat A$ is a \emph{lift candidate} for $(A, \{e,f,g\})$. Recall that $A-U$ denotes the matrix obtained from $A$ by removing the rows and columns labelled by elements of $U$.
  \begin{aclaim}\label{cl:threeleavespivot}
    If $(A,\{e,f,g\})$ is a bad pair with lift candidate $\hat A$, and $x\in X$, $y\in Y$ are such that $A_{xy}\neq 0$ and $\{x,y\}\cap\{e,f,g\}=\emptyset$, then $(A^{xy},\{e,f,g\})$ is a bad pair with lift candidate $\hat A^{xy}$.
  \end{aclaim}
  \begin{subproof}
    Since $A^{xy}$ has a global lift if and only if $A$ has a global lift, $A^{xy}$ is a minimal counterexample to the theorem. Since $\bip(A-U)$ is connected for all $U\subseteq \{e,f,g\}$, Lemma~\ref{lem:bipconn}\eqref{bip:conn} implies that $\bip(A^{xy}-U)$ is connected for all $U\subseteq \{e,f,g\}$. A spanning tree $T'$ for $A^{xy}$ with leaves $\{e,f,g\}$ is now easily found, so $(A,\{e,f,g\})$ is indeed a bad pair. Pivoting commutes with deleting rows and columns other than $x$, $y$. From this and the fact that $\hat A - U$ is a global $^\uparrow$-lift of $A-U$ for all nonempty $U \subseteq\{e,f,g\}$ it follows that $\hat A^{xy}$ is a lift candidate for $(A^{xy},\{e,f,g\})$.
  \end{subproof}
  We say that $(A,\{e,f,g\})$ is a \emph{local} bad pair if a lift candidate $\hat A$ is not a local lift of $A$.
  In that case there exist $X'\subseteq X$, $Y'\subseteq Y$, $|X'| = |Y'|$, such that either
  \begin{enumerate}
    \item $\det(\hat A[X',Y'])$ is undefined, or
    \item $\bip(A[X',Y'])$ is a cycle $C$ but $\sigma_{\hat A}(C) \neq \sigma_A(C)^\uparrow$.
  \end{enumerate}
  We call $(X',Y')$ a \emph{certificate}.
  \begin{aclaim}
    If there exists a counterexample $A$ to the theorem with $|X|+|Y|$ minimal such that $A$ has no local lift then there exist $e,f,g \in X\cup Y$ such that one of $(A, \{e,f,g\})$ and $(A^T, \{e,f,g\})$ is a bad pair.
  \end{aclaim}
  \begin{subproof}
    Let $A$ be a counterexample to the theorem with $|X|+|Y|$ minimal such that $A$ has no local lift. By Claim~\ref{cl:lifttwoconn} $\bip(A)$ is 2-connected.
    From Lemma~\ref{lem:signature}\eqref{sgn:det} it follows that $\bip(A)$ is not a cycle. By Lemma~\ref{lem:biptree} there exists a spanning tree $T$ of $\bip(A)$ which has leaves $e,f,g$, with $e,f \in X$ and $g \in Y$ or $e,f \in Y$ and $g \in X$. Clearly if $A$ is a counterexample then so is $A^T$. The claim follows.
  \end{subproof}
  \begin{aclaim}
    Let $(A,\{e,f,g\})$ be a local bad pair with certificate $(X',Y')$ such that $|X'|$ is minimal. Then $|X'| = 2$ and all entries of $A[X',Y']$ are nonzero.
  \end{aclaim}
  \begin{subproof}
    By Claim~\ref{cl:minorminimal} we have $X'\cup Y' \supseteq \{e,f,g\}$ so $|X'| \geq 2$. If there are $x \in X'\setminus\{e,f\}$, $y \in Y'\setminus g$ with $A_{xy} \neq 0$ then it follows from Claim~\ref{cl:threeleavespivot} and one of Lemma~\ref{lem:detpivot} and Lemma~\ref{lem:signature}\eqref{sgn:pivot} that $(A^{xy}, \{e,f,g\})$ is a bad pair with lift candidate $\hat A^{xy}$ and certificate $(X'\setminus x, Y'\setminus y)$, which contradicts the minimality of $|X'|$.

    If there is an $x \in X'\setminus\{e,f\}$ then $A_{xy} = 0$ for all $y \in Y'\setminus \{g\}$. Then $\det(\hat A[X',Y']) = \hat A_{xg}\det(\hat A[X'\setminus x, Y'\setminus g])$. But $\hat A-\{x,g\}$ is a square submatrix of $\hat A-g$ so its determinant is defined, a contradiction. It follows that $|X'| = |Y'| = 2$.

    If some entry of $\hat A[X',Y']$ equals 0 then clearly $\bip(A[X',Y'])$ is not a cycle, so $\det(\hat A[X',Y'])$ must be undefined. But this determinant is the product of entries in $\hat A$ and, possibly, $-1$. This is a contradiction since all entries are in $\hat\parf$. The claim follows.
  \end{subproof}
  Suppose $(A,\{e,f,g\})$ is a local bad pair with minimal certificate $(X',Y')$, i.e. $|X'| = 2$. Suppose $X' =\{e,f\}, Y' = \{g,h\}$. Since all four entries of $\hat A[X',Y']$ are nonzero, clearly $\sigma_{\hat A}(C) \neq \sigma_A(C)^\uparrow$ for $C = (e,g,f,h,e)$.
  \begin{aclaim}
    If $(A, \{e,f,g\})$ is a local bad pair with minimal certificate then there exist $p,q,r,s \in \parf$ such that $A$ is scaling-equivalent to one of the following matrices:
    \begin{align}\label{eq:badmatrices}
      A_1 :=  \kbordermatrix{& h & g\\
                          i &  \textcircled{1} & \textcircled{1}\\
                          e & \textcircled{1} & p\\
                          f &  \textcircled{1} & q
             },\quad
      A_2 :=  \kbordermatrix{& j & h & g\\
                                 i & \textcircled{1} & 0 & \textcircled{1}\\
                                 k & \textcircled{1} & \textcircled{1} & 0\\
                                 e & p & \textcircled{1} & r\\
                                 f & q & \textcircled{1} & s
             }.
    \end{align}
  \end{aclaim}
  \begin{subproof}
    Let $(X',Y')$ be a minimal certificate, say $X' =\{e,f\}$ and $Y' = \{g,h\}$ for some $g \in Y$.
    Since $\bip(A - \{e,f\})$ is connected, there exists a $g-h$ path $P$ in $\bip(A -\{e,f\})$. Let $P$ be a shortest such path. Then $\bip(A[V(P)]) = P$. Then $T := P \cup \{he,hf\}$ is a spanning tree for $A' := A[V(P)\cup\{e,f\}]$ with leaves $\{e,f,g\}$. But if $\hat A'$ is a lift candidate for $(A',\{e,f,g\})$, then $\hat A[V(P)\cup\{e,f\}] \sim \hat A'$ by Lemma~\ref{lem:uniqA2}, so $(A', \{e,f,g\})$ is a local bad pair with certificate $(\{e,f\},\{g,h\})$. By the minimality of $|X|+|Y|$ we then have $A = A'$.

    If $|V(P)|\geq 7$ then $P$ has an edge $xy$ with $x\in X$ such that $A_{xg} = A_{xh} = 0$. By Claim~\ref{cl:threeleavespivot} we have that $(A^{xy}, \{e,f,g\})$ is a local bad pair with minimal certificate. But $A^{xy}$ has a shorter $g-h$ path, which again contradicts the minimality of $|X|+|Y|$. Therefore $|V(P)| = 3$ or $|V(P)| = 5$, from which the claim follows.
  \end{subproof}
  \begin{aclaim}\label{cl:locallift}
    There does not exist a local bad pair.
  \end{aclaim}
  \begin{subproof}
    Suppose $(A, \{e,f,g\})$ is a local bad pair with minimal certificate. Since \eqref{enum:propfour} does not hold we have $A\not\sim A_1$. Therefore $A \sim A_2$. Assume, without loss of generality, that $A = A_2$ for some $p,q,r,s$. Let $\hat p, \hat q, \hat r, \hat s$ be the entries of $\hat A$ corresponding to $p,q,r,s$.
    \begin{sclaim}\label{scl:notbothzero}$p$ and $q$ are not both zero.\end{sclaim}
    \begin{subproof}
      $A^{ij} - \{i,j\}$ is scaling-equivalent to a matrix of the form $A_1$, a contradiction.
    \end{subproof}
    \begin{sclaim}\label{scl:notbothnonzero}
      Either $p = 0$ or $q = 0$.
    \end{sclaim}
    \begin{subproof}
      Suppose $p\neq 0$, $q \neq 0$. Then $\hat p = p^\uparrow, \hat q = q^\uparrow, \hat r = (r/p)^\uparrow p^\uparrow$, and $\hat s = (s/q)^\uparrow q^\uparrow$. Since $\sigma_{\hat A}(C) \neq \sigma_A(C)^\uparrow$ for $C = (e,g,f,h,e)$ it follows that
    \begin{align}\label{eq:cyclewrong}
      \frac{\hat r}{\hat s} \neq \left(\frac{r}{s}\right)^\uparrow.
    \end{align}
    $A$ is minor-minimal, so $A[\{e,f\},\{j,h,g\}]$ has a local $^\uparrow$-lift. This matrix is scaling-equivalent to the following normalized matrices:
    \begin{align}
      \kbordermatrix{
          & j & h & g\\
        e & \textcircled{1} & \textcircled{1} & r/s\\
        f & q/p & \textcircled{1} & \textcircled{1}
      }, \quad
      \kbordermatrix{
          & j & h & g\\
        e & \textcircled{1} & \textcircled{1} & \textcircled{1}\\
        f & \textcircled{1} & p/q & \frac{ps}{qr}
      }.
    \end{align}
    Since these matrices have a local $^\uparrow$-lift we conclude, using $(1/p)^\uparrow = 1/(p^\uparrow)$, that
    \begin{align}
      \left(\frac{p}{q}\right)^\uparrow\left(\frac{s}{r}\right)^\uparrow = \left(\frac{ps}{qr}\right)^\uparrow.
    \end{align}
    Likewise $A[\{i,e,f\},\{j,g\}]$ has a local $^\uparrow$-lift. This gives
    \begin{align}
      \left(\frac{p}{r}\right)^\uparrow\left(\frac{s}{q}\right)^\uparrow = \left(\frac{ps}{qr}\right)^\uparrow.
    \end{align}
    Finally, $A_1[\{k,e,f\},\{j,h\}]$ has a local $^\uparrow$-lift. This gives
    \begin{align}
      \frac{p^\uparrow}{q^\uparrow} = \left(\frac{p}{q}\right)^\uparrow.
    \end{align}
    But then
    \begin{align}
      \left(\frac{r}{s}\right)^\uparrow = \left(\frac{r}{p}\right)^\uparrow p^\uparrow / \left( \left(\frac{s}{q}\right)^\uparrow q^\uparrow \right) = \frac{\hat r}{\hat s},
    \end{align}
    a contradiction.
    \end{subproof}
    By symmetry we may assume $p = 0$.
    \begin{sclaim}\label{scl:qone}
      $q = 1$.
    \end{sclaim}
    \begin{subproof}
      Suppose $p = 0$, $q \neq 0$, $q \neq 1$. Then $A^{kh}$ is scaling-equivalent to
      \begin{align}
        A' := \kbordermatrix{& j & k & g\\
                                 i & \textcircled{1} & 0 & \textcircled{1}\\
                                 h & \textcircled{1} & \textcircled{1} & 0\\
                                 e & p' & \textcircled{1} & r'\\
                                 f & q' & \textcircled{1} & s'
              }
      \end{align}
      with $p' = 1$, $q' = 1-q$, $r'=-r$, $s'=-s$. A spanning tree $T'$ has been circled. Let $\hat A'$ be a $T$-normalized lift candidate for $(A,\{e,f,g\})$. By Claim~\ref{cl:threeleavespivot} $\hat A'\sim\hat A^{kh}$.
      But $\hat A'[\{e,f\},\{k,g\}]$ is, after exchanging the labels $k$ and $h$, scaling-equivalent to $\hat A[\{e,f\},\{h,g\}]$, so  again $\sigma_{\hat A'}(C) \neq \sigma_{A'}(C)^\uparrow$ for $C = (e,g,f,k,e)$, by Lemma~\ref{lem:signature}\eqref{sgn:scale}. But this is impossible by Claim~\ref{scl:notbothnonzero}.
    \end{subproof}
    Now $p = 0$, $q = 1$. Then $\hat s = s^\uparrow$ and $\hat r = -(-r)^\uparrow$. Scale row $e$ of $A$ by $1/r$ and then column $h$ by $r$. After permuting some rows and columns we obtain
    \begin{align}
      A' :=                \kbordermatrix{& g & j & h\\
                                   e & \textcircled{1} & 0 & \textcircled{1}\\
                                   i & \textcircled{1} & \textcircled{1} & 0\\
                                   k & 0 & \textcircled{1} & r\\
                                   f & s & \textcircled{1} & r
             }.
    \end{align}
    A spanning tree $T'$ has been circled. Let $\hat A'$ be the $T'$-normalized lift candidate for $(A',\{k,f,h\})$. Then $\hat A'_{kh} = -(-r)^\uparrow$ and $\hat A'_{fh} = (r/s)^\uparrow s^\uparrow$. But then $\sigma_{\hat A}(C) \neq \sigma_A(C)^\uparrow$ for $C' = (k,j,f,h,k)$. By Claim~\ref{scl:qone} we have $s = 1$. We can now repeat the argument and conclude that also $r = 1$. Hence \eqref{enum:propfour} holds, contradicting our choice of $A$.
    This ends the proof of Claim~\ref{cl:locallift}.
  \end{subproof}
  A pair $(A, xy)$, where $A$ is an $X\times Y$ $\parf$-matrix such that $X\cap Y = \emptyset$, and $x \in X$, $y \in Y$ are such that $A_{xy} \neq 0$, is called a \emph{bad-pivot pair} if
  \begin{enumerate}
    \item $A$ is a counterexample to the theorem with $|X|+|Y|$ minimal;
    \item $A$ has a local lift $\hat A$, but $\hat A^{xy}$ is not a local lift of $A^{xy}$.
  \end{enumerate}
  \begin{aclaim}
    There exists a bad-pivot pair.
  \end{aclaim}
  \begin{subproof}
    Let $A$ be a counterexample to the theorem with $|X|+|Y|$ minimal. By Claim~\ref{cl:locallift} $A$ has a local lift $\hat A$. Suppose $\hat A$ is not a global $^\uparrow$-lift for $A$. Then there exist sequences $A_0, \ldots, A_k$ and $\hat A_0,\ldots, \hat A_k$ such that $A_0 = A$, $\hat A_0 = \hat A$, and for $i = 1, \ldots, k$, $A_i = (A_{i-1})^{x_iy_i}$ and $\hat A_i = (\hat A_{i-1})^{x_iy_i}$, so that $\hat A_k$ is not a local $^\uparrow$-lift of $A_k$. Choose $A$ and these sequences such that $k$ is as small as possible. But then $k=1$, so there is an edge $xy\in\bip(A)$ such that $A_{xy} \neq 0$ and $\hat A^{xy}$ is not a local $^\uparrow$-lift of $A^{xy}$.
  \end{subproof}
  By Claim~\ref{cl:threeleavespivot} we have
  \begin{aclaim}
    If $(A, \{e,f,g\})$ is a bad pair and $(A, xy)$ is a bad-pivot pair, then $\{x,y\}\cap\{e,f,g\}\neq \emptyset$.
  \end{aclaim}
  Let $T'$ be a tree such that $x, y \in T'$ and $T'$ has three leaves $\{e',f',g'\}$, not all rows and not all columns, such that $\{x,y\}\cap\{e',f',g'\} = \emptyset$. From the proof of Lemma~\ref{lem:biptree} we conclude that we can extend $T'$ to a spanning tree of $\bip(A)$ with three leaves $\{e,f,g\}$, not all rows and not all columns, such that $\{x,y\}\cap\{e,f,g\} = \emptyset$. We call $T'$ ``good for $xy$''. It follows that there is no good tree for $xy$ in $\bip(A)$.
  \begin{aclaim}
    There exists a bad-pivot pair $(A, xy)$ such that, for some $p, q \in \parf$, we have
    \begin{align}\label{eq:badpivotpair}
      A = \kbordermatrix{
         & y & g & h\\
       x & \textcircled{1} & \textcircled{1} & 0\\
       e & \textcircled{1} & p & \textcircled{1}\\
       f & 0 & \textcircled{1} & q
      }.
    \end{align}
  \end{aclaim}
  \begin{subproof}
    Let $(A,xy)$ be a bad-pivot pair.
    By Claim~\ref{cl:lifttwoconn} $\bip(A)$ is 2-connected, so there exists a cycle $C$ containing $xy$. 
    By Lemma~\ref{lem:signature}\eqref{sgn:pivot},\eqref{sgn:det} $\bip(A)$ is not a cycle. Then there exists a path $P$ between two vertices of $C$, which is internally vertex-disjoint from $C$.
    If some vertex $v\in P\cap C$ is not in $\delta(\{x,y\})$ then we delete the two edges of $C$ adjacent to $v$ and obtain a good tree for $xy$, a contradiction. If $x \in P \cap C$ then we delete an edge of $C$ not adjacent to $xy$ and an edge of $P$ not adjacent to $xy$ to obtain a good tree for $xy$, a contradiction. Since $\bip(A)$ is simple and bipartite, such edges exist. Therefore we may assume that all such paths $P$ have the neighbours $u$, $v$ of $xy$ as end vertices. If $P$ has length at least 3 and $C$ has length at least 6 then again a good tree for $xy$ can be found. If  $P$ has length at least 3 and $C$ has length $4$, then we can replace $C$ by $C' := C\setminus uv \cup P$, and $P$ by $P' := uv$. Therefore, without loss of generality, we may assume $P$ has length 1.

    Assume a bad-pivot pair $(A, xy)$ was chosen such that the length of $P$ is 1 and the length of $C$ is as small as possible. Suppose $C$ has length more than 6. Let $x'y'$ be the edge of $C$ at maximum distance from $xy$. We can find a good tree for $x'y'$, so $\hat A' := \hat A^{x'y'}$ is a local $^\uparrow$-lift of $A' := A^{x'y'}$. But in $\bip(A')$ there is a good tree for $xy$, so $(\hat A')^{xy}$ is a local lift for $(A')^{xy}$. But $((\hat A')^{xy})^{y'x'} = \hat A^{xy}$, so there is no good tree for $y'x'$ in $(A')^{xy}$. This is only the case if $A^{xy}$ is a cycle. But it is easily checked that in this case $\hat{A^{xy}} = \hat A^{xy}$, a contradiction. The claim follows.
  \end{subproof}
    Suppose $(A, xy)$ is a bad-pivot pair with $A$ as in \eqref{eq:badpivotpair} for some $p,q\in\parf$. The normalized local $^\uparrow$-lift $\hat A$ of $A$ has $\hat A_{eg} = p^\uparrow$ and $\hat A_{fh} = (pq)^\uparrow/ p^\uparrow$. After a pivot over $xy$ and renormalization we have
    \begin{align}
      A' = \kbordermatrix{
         & x & g & h\\
       y & \textcircled{1} & \textcircled{1} & 0\\
       e & \textcircled{1} & 1-p & \textcircled{1}\\
       f & 0 & \textcircled{1} & -q
      }.
    \end{align}
    The normalized local $^\uparrow$-lift $\hat A'$ of $A'$ has $\hat A'_{eg} = (1-p)^\uparrow$ and $\hat A'_{fh} = (q(p-1))^\uparrow/ (1-p)^\uparrow$.
    By definition of the lifting function $(1-p)^\uparrow = 1 - p^\uparrow$ and $\left(\frac{p}{p-1}\right)^\uparrow = \frac{p^\uparrow}{p^\uparrow -1}$. Since $\hat A'$ is not scaling-equivalent to $\hat A^{xy}$, we must have
    \begin{align}\label{eq:globliftproblem}
      -(pq)^\uparrow/ p^\uparrow \neq (q(p-1))^\uparrow/(1-p)^\uparrow.
    \end{align}
    Consider
    \begin{align}
      A^{xg} = \kbordermatrix{
         & y & x & h\\
       g & 1 & -1 & 0\\
       e & 1-p & p & 1\\
       f & -1 & 1 & q
      }.
    \end{align}
    Since $A$ is minor-minimal, $A^{xg}[\{e,f\},\{y,x,h\}]$ has a global $^\uparrow$-lift. If we normalize with respect to tree $T' = \{ey,ex,eh,fy\}$ then we find
    \begin{align}
      \left(\frac{p-1}{p}\right)^\uparrow(pq)^\uparrow = ((1-p)q)^\uparrow
    \end{align}
    which contradicts \eqref{eq:globliftproblem}. Therefore $A$ does have a global $^\uparrow$-lift.
  It follows that no counterexample exists, which completes the proof of the theorem.
\end{proof}
\addtocounter{theorem}{2}
We remark here that for most of our applications, including all examples in the next section,  the restriction of $\phi$ to the fundamental elements, denoted $\phi|_{\fun(\hat\parf)}$, is a bijection between $\fun(\hat\parf)$ and $\fun(\parf)$. Then $(\phi|_{\fun(\hat\parf)})^{-1}$ is an obvious choice for the lifting function.
We did not specify this lifting function in the theorem statement because we need the more general version for the proof of Lemma~\ref{lem:liftpf}.

We have the following corollary:
\begin{corollary}\label{cor:pfequivalence}
  Let $\parf$, $\hat\parf$, $\phi$, $^\uparrow$ be as in Theorem~\ref{thm:lift}. Suppose that
  \begin{enumerate}
    \item\label{equi:one} If $1+1 \doteq 0$ in $\parf$ then $1+1 \doteq 0$ in $\hat \parf$;
    \item\label{equi:two} If $1+1$ is defined and nonzero in $\parf$ then $1+1$ is defined and nonzero in $\hat \parf$;
    \item\label{equi:triple} For all $p,q,r \in \fun(\parf)$ such that $pqr = 1$, we have $p^\uparrow q^\uparrow r^\uparrow = 1$.
  \end{enumerate}
  Then a matroid is $\parf$-representable if and only if it is $\hat\parf$-representable.
\end{corollary}
\begin{proof}
Since there is a nontrivial homomorphism $\phi:\hat\parf\rightarrow\parf$, every matroid that is $\hat\parf$-representable is also $\parf$-representable. To prove the other implication it suffices to show that every $\parf$-matrix has a global $^\uparrow$-lift. Suppose that this is false. By Theorem~\ref{thm:lift} there must be a $\parf$-matrix $B$ as in \eqref{eq:mincounterexmats} that does not have a local $^\uparrow$-lift. Suppose there are $p',q' \in \parf$ such that the following $\parf$-matrix has no local $^\uparrow$-lift:
\begin{align}\label{eq:badmatrixone}
  \begin{bmatrix}
    1 & 1 & 1\\
    1 & p' & q'
  \end{bmatrix}.
\end{align}
This matrix has a local $^\uparrow$-lift if and only if
\begin{align}\label{eq:liftUtwofive}
  \left(\frac{p'}{q'}\right)^\uparrow = \frac{(p')^\uparrow}{(q')^\uparrow}.
\end{align}
Pick $p := p'$, $q := (q')^{-1}$, and $r := q'/p'$. Then \eqref{eq:liftUtwofive} holds if and only if $p^\uparrow q^\uparrow r^\uparrow = 1$, which follows from \eqref{equi:triple}. It follows that
\begin{align}
  A = \begin{bmatrix}
    0 & 1 & 1 & 1\\
    1 & 0 & 1 & 1\\
    1 & 1 & 0 & 1
  \end{bmatrix}
\end{align}
has no local $^\uparrow$-lift. Note that $A$ has a cycle having signature $-1$. Hence $1-(-1)^\uparrow$ must be defined in $\parf$, and hence also in $\hat\parf$. Since $\phi(1) +\phi((-1)^\uparrow) \doteq 0$, we have $(-1)^\uparrow = -1$. Moreover, \eqref{equi:one} and \eqref{equi:two} imply that $1+1 \doteq 0$ in $\parf$ if and only if $1+1 \doteq 0$ in $\hat\parf$, since $\phi(1) = 1$. Let $\hat A$ be a $\hat\parf$-matrix such that $\hat A_{xy} = 1$ if $A_{xy} = 1$ and $\hat A_{xy} = 0$ if $A_{xy} = 0$. It is easily checked that all conditions of Definition~\ref{def:locallift} are met, so $\hat A$ is a local $^\uparrow$-lift of $A$, a contradiction.
\end{proof}

\section{Applications}\label{sec:examples}
In this section we use the notation related to fundamental elements that was introduced in Section~\ref{ssec:funds}.

\subsection{Binary matroids}
In addition to Theorem~\ref{thm:reg}, Tutte~\cite{Tut65} proved the following characterization of regular matroids:

\begin{theorem}\label{thm:regex}
  Let $M$ be a binary matroid. Exactly one of the following is true:
  \begin{enumerate}
    \item $M$ is regular;
    \item $M$ has a minor isomorphic to one of $F_7$ and $F_7^*$.
  \end{enumerate}
\end{theorem}

The shortest known proof for this result is by Gerards~\cite{Ger89}. The techniques used to prove the lift theorem generalize those used by Gerards, so it is no surprise that Theorem~\ref{thm:regex} can also be proven using the Lift Theorem. Recall from Definition~\ref{prop:pffromring} that $\parf(\ring, S)$ is the partial field $(\langle S\cup \{-1\}\rangle, +, \cdot, 0, 1)$, where multiplication and addition are the restriction of the operations in $\ring$.

\begin{proof}
  Let $\parf := \GF(2)$, let $\hat\parf := \reg = \parf(\Q, \{-1,0,1\})$, let $\phi:\hat\parf\rightarrow\parf$ be defined by $\phi(-1) = \phi(1) = 1$, $\phi(0) = 0$, and let $^\uparrow:\fun(\parf)\rightarrow\hat\parf$ be defined by $0^\uparrow = 0$, $1^\uparrow = 1$. It is readily checked that this is a lifting function.

  It is not hard to see that $F_7$ and $F_7^*$ are not regular. For the converse, let $M$ be a binary matroid without $F_7$- and $F_7^*$-minor, and let $A$ be a $\parf$-matrix such that $M = M[I|A]$. All rank-2 binary matroids are regular, so $A$ has no minor isomorphic to a matrix as in \eqref{eq:mincounterexmats}. But then Theorem~\ref{thm:lift} implies that $A$ has a global $\hat\parf$-lift, and hence $M$ is regular.
\end{proof}

Tutte proved Theorem~\ref{thm:regex} using his Homotopy Theorem~\cite{Tuthom}. We believe that the Homotopy Theorem can be used to prove the Lift Theorem as well.

\subsection{Ternary matroids}
Our first applications of the Lift Theorem consist of new proofs of three results of Whittle~\cite{Whi97}.

First we prove Theorem~\ref{thm:dyadicintro} from the introduction. A matroid is called \emph{dyadic} if it is representable over the partial field $\dyadic := \parf(\Q,2)$. First we compute the set of fundamental elements. Recall that $\assoc\{0\} = \assoc\{1\} = \{0,1\}$, and $\assoc\{p\} =\left\{p,1-p,\frac{1}{1-p},\frac{p}{p-1},\frac{p-1}{p},\frac{1}{p}\right\}$.

\begin{lemma}$\fun(\dyadic) = \assoc\{1,2\} = \{0,1,-1,2,1/2\}$.
\end{lemma}
\begin{proof}
  We find all solutions of
  \begin{align}
    1-p = q
  \end{align}
  where $p = (-1)^s 2^x$ and $q = (-1)^t 2^y$. If $x < 0$ then we divide both sides by $p$. Likewise if $y<0$ then we divide both sides by $q$. We may multiply both sides with $-1$. After rearranging and dividing out common factors we need to find all solutions of
  \begin{align}
    2^{x'} + (-1)^{s'} 2^{y'} + (-1)^{t'} = 0
  \end{align}
  where $x',y' \geq 0$. This equation has solutions only if one of $2^{x'}$, $2^{y'}$ is odd. This implies that we just need to find all solutions of
  \begin{align}
    2^{x''} +(-1)^{s''} + (-1)^{t''} = 0.
  \end{align}
  There are finitely many solutions. Enumeration of these completes the proof.
\end{proof}
\begin{theorem}[Whittle~\cite{Whi97}]\label{thm:dyadic}
Let $M$ be a matroid. The following are equivalent:
\begin{enumerate}
  \item \label{eq:Dthreefive}$M$ is representable over $\GF(3)\otimes\GF(5)$;
  \item \label{eq:DTU} $M$ is $\dyadic$-representable;
  \item \label{eq:Devery} $M$ is representable over every field that does not have characteristic 2.
\end{enumerate}
\end{theorem}
\begin{proof}Let $\phi_3:\dyadic\rightarrow\GF(3)$ be determined by $\phi(2) = -1$. Let $\phi_5:\dyadic\rightarrow\GF(5)$ be determined by $\phi(2) = 2$. Clearly both are partial field homomorphisms. But then $\phi = \phi_3\otimes\phi_5$ is a partial field homomorphism $\dyadic\rightarrow\GF(3)\otimes\GF(5)$. It is readily seen that $\phi|_{\fun(\dyadic)}:\fun(\dyadic)\rightarrow\fun(\GF(3)\otimes\GF(5))$ is bijection. Taking $(\phi|_{\fun(\dyadic)})^{-1}$ as lifting function we apply Corollary \ref{cor:pfequivalence}, thereby proving \eqref{eq:Dthreefive}$\Leftrightarrow$\eqref{eq:DTU}. For \eqref{eq:DTU}$\Rightarrow$\eqref{eq:Devery}, use again suitable homomorphisms. The implication \eqref{eq:Devery}$\Rightarrow$\eqref{eq:Dthreefive} is trivial, by Corollary~\ref{thm:listoffields}.
\end{proof}

A matroid is called \emph{near-regular} if it is representable over the partial field $\nreg := \parf(\Q(\alpha),\{\alpha,1-\alpha\})$, where $\alpha$ is an indeterminate.

\begin{lemma}$\fun(\nreg) = \assoc\{1,\alpha\}$.
\end{lemma}
\begin{proof}
  We find all $p = (-1)^s \alpha^x (1-\alpha)^y$ such that $1-p\doteq q$ in $\nreg$. Consider the homomorphism $\phi:\nreg \rightarrow \dyadic$ determined by $\phi(\alpha) = 2$. Since fundamental elements must map to fundamental elements, it follows that $x \in \{-1,0,1\}$. Likewise, $\psi: \nreg \rightarrow \dyadic$, determined by $\psi(\alpha)=-1$, shows that $y \in \{-1,0,1\}$. Again, a finite check remains.
\end{proof}

\begin{theorem}[Whittle~\cite{Whi97}]\label{thm:nreg}
Let $M$ be a matroid. The following are equivalent:
\begin{enumerate}
  \item \label{eq:nregthreefourfive}$M$ is representable over $\GF(3)\otimes\GF(4)\otimes\GF(5)$;
  \item \label{eq:nregthreeeight}$M$ is representable over $\GF(3)\otimes\GF(8)$;
  \item \label{eq:nregTU} $M$ is $\nreg$-representable;
  \item \label{eq:nregevery} $M$ is representable over every field with at least 3 elements.
\end{enumerate}
\end{theorem}
\begin{proof}Let $\phi:\nreg\rightarrow \GF(3)\otimes\GF(4)\otimes\GF(5)$ be determined by $\phi(\alpha)=(-1,\omega,2)$. Again $\phi|_{\fun(\nreg)}:\fun(\nreg)\rightarrow\fun(\GF(3)\otimes\GF(4)\otimes\GF(5))$ is a bijection, so we use $(\phi|_{\fun(\nreg)})^{-1}$ as lifting function and apply Corollary \ref{cor:pfequivalence} to prove \eqref{eq:nregthreefourfive}$\Leftrightarrow$\eqref{eq:nregTU}. For  \eqref{eq:nregTU}$\Rightarrow$\eqref{eq:nregevery}, use a homomorphism $\phi'$ such that $\phi'(\alpha) = p$ for any $p \in \field\setminus\{ 0, 1\}$. Similar constructions prove the remaining implications.
\end{proof}

Let $\splittable := \parf(\C,\{2,\zeta\})$, where $\zeta$ is a primitive complex sixth root of unity.
\begin{lemma}$\fun(\splittable) = \assoc\{1,2,\zeta\} = \{0,1,-1,2,1/2,\zeta,1-\zeta\}$.
\end{lemma}
\begin{proof}
  Clearly all these elements are fundamental elements. The complex argument of every element of $\splittable$ is equal to a multiple of $\pi/3$, from which it follows easily that no other fundamental elements exist.
\end{proof}

\begin{theorem}[Whittle~\cite{Whi97}]\label{thm:splittable}
Let $M$ be a matroid. The following are equivalent:
\begin{enumerate}
  \item \label{eq:splittablethreeseven}$M$ is representable over $\GF(3)\otimes\GF(7)$;
  \item \label{eq:splittableTU} $M$ is $\splittable$-representable;
  \item \label{eq:splittableCong} $M$ is representable over $\GF(3)$, over $\GF(p^2)$ for all primes $p > 2$, and over $\GF(p)$ when $p \equiv 1 \mod 3$.
\end{enumerate}
\end{theorem}
\begin{proof}Let $\phi:\splittable\rightarrow \GF(3)\otimes\GF(7)$ be determined by $\phi(2) = (-1,2)$ and $\phi(\zeta) = (-1,3)$. Again $\phi|_{\fun(\splittable)}:\fun(\splittable)\rightarrow\fun(\GF(3)\otimes\GF(7))$ is a bijection, so we use $(\phi|_{\fun(\splittable)})^{-1}$ as lifting function and apply Corollary \ref{cor:pfequivalence} to prove \eqref{eq:splittablethreeseven}$\Leftrightarrow$\eqref{eq:splittableTU}.   For \eqref{eq:splittableTU}$\Rightarrow$\eqref{eq:splittableCong} we use an argument similar to the proof of Theorem~\ref{thm:sru}. Note that the ring $\Z[\frac{1}{2}, \zeta]$ is not the ring of integers of an algebraic number field, but every element is of the form $2^k x$ for some $k \in \Z$, $x \in \Z[\zeta]$. Hence, in contrast to the partial field $\psru$, there are no homomorphisms to finite fields of characteristic 2. Finally, \eqref{eq:splittablethreeseven} is a special case of \eqref{eq:splittableCong}.
\end{proof}

\subsection{Quaternary and quinary matroids}
Our next example is a proof of Theorem~\ref{thm:golratintro}. A matroid is called \emph{golden ratio} (in \cite{Whi05} ``golden mean'' is used) if it is representable over the partial field $\golrat := \parf(\R, \tau)$, where $\tau$ is the golden ratio, i.e. the positive root of $x^2-x-1=0$.

\begin{lemma}$\fun(\golrat) = \assoc\{1,\tau\} = \{0,1,\tau, -\tau, 1/\tau, -1/\tau, \tau^2, 1/\tau^2\}$.
\end{lemma}
\begin{proof}
  Remark that for all $k \in \Z$, $\tau^k = f_{k} + f_{k+1} \tau$, where $f_0 = 0$, $f_1 = 1$, and $f_{i+2} - f_{i+1} - f_i = 0$, i.e. the Fibonacci sequence, extended to hold for negative $k$ as well. If $p = (-1)^s(f_k + f_{k+1}\tau)$ is a fundamental element, then $\{|(-1)^sf_{k}-1|,|f_{k+1}|\}$ has to be a set of two consecutive Fibonacci numbers. We leave out the remaining details.
\end{proof}

\begin{theorem}[Vertigan]\label{thm:golrat}
Let $M$ be a matroid. The following are equivalent:
\begin{enumerate}
  \item \label{eq:GRfourfive}$M$ is representable over $\GF(4)\otimes\GF(5)$;
  \item \label{eq:GRTU} $M$ is $\golrat$-representable;
  \item \label{eq:GRevery} $M$ is representable over $\GF(5)$, over $\GF(p^2)$ for all primes $p$, and over $\GF(p)$ when  $p \equiv \pm 1 \mod 5$.
\end{enumerate}
\end{theorem}
\begin{proof}Let $\phi:\golrat\rightarrow \GF(4)\otimes\GF(5)$ be determined by $\phi(\tau)= (\omega,3)$. Again $\phi|_{\fun(\golrat)}:\fun(\golrat)\rightarrow\fun(\GF(4)\otimes\GF(5))$ is a bijection, so we use $(\phi|_{\fun(\golrat)})^{-1}$ as lifting function and apply Corollary \ref{cor:pfequivalence} to prove \eqref{eq:GRfourfive}$\Leftrightarrow$\eqref{eq:GRTU}.

For \eqref{eq:GRTU}$\Rightarrow$\eqref{eq:GRevery} we use an argument similar to the proof of Theorem~\ref{thm:sru}. Finally, \eqref{eq:GRfourfive} is a special case of \eqref{eq:GRevery}.
\end{proof}

A matroid is called \emph{Gaussian} if it is representable over the partial field $\gauss := \parf(\C,\{i,1-i\})$, where $i$ is a root of $x^2+1=0$.

\begin{lemma}\label{lem:gaussfun}
\begin{align}
  \fun(\gauss) = \assoc\{1,2,i\} = \left\{0,1,-1,2,\tfrac{1}{2},i,i+1,\tfrac{i+1}{2},1-i,\tfrac{1-i}{2},-i\right\}.
\end{align}
\end{lemma}
\begin{proof}
  First note that the complex argument of every element of $\gauss$ is a multiple of $\pi/4$. It follows that if $p = i^x (1-i)^y$ is a fundamental element, then $\frac{1}{\sqrt{2}} \leq p \leq \sqrt{2}$.  Therefore there are finitely many fundamental elements in $\C\setminus\R$. It is easily checked that all numbers on the real line are powers of 2. The result follows.
\end{proof}
Our next result requires more advanced techniques. The following lemma is a corollary of Whittle's Stabilizer Theorem~\cite{Whi96b}.
\begin{theorem}[Whittle~\cite{Whi96b}]\label{Utwofivestabilizer}Let $M$ be a 3-connected quinary matroid with a minor $N$ isomorphic to one of $U_{2,5}$ and $U_{3,5}$. Then any representation of $M$ over $\GF(5)$ is determined up to strong equivalence by the induced representation of $N$.
\end{theorem}
\begin{lemma}\label{lem:hydratwo}Let $M$ be a 3-connected matroid.
\begin{enumerate}
  \item\label{eq:onehydratwo} If $M$ has at least $2$ inequivalent representations over $\GF(5)$, then $M$ is representable over $\hydra_2$.
  \item\label{eq:twohydratwo} If $M$ has a $U_{2,5}$- or $U_{3,5}$-minor and $M$ is representable over $\hydra_2$, then $M$ has at least $2$ inequivalent representations over $\GF(5)$.
\end{enumerate}
\end{lemma}
\begin{proof}
  Let $\phi:\gauss\rightarrow \GF(5)\otimes\GF(5)$ be determined by $\phi(i) = (2,3)$. Then $\phi(2) = \phi(i(1-i)^2) = (2,2)$. Let $\phi_i:\GF(5)\otimes\GF(5)\rightarrow\GF(5)$ be determined by $\phi_i(x) = x_i$ for $i = 1,2$. Let
  \begin{align}
    A:=\begin{bmatrix}1 & 1 & 1\\1 & p' & q'\end{bmatrix}
  \end{align}
  for some, $p',q' \in \gauss$. If $A$ is an $\gauss$-matrix then $p',q'\in\fun(\gauss)$. A finite check then shows that for each of these,  $\phi_1(\phi(A))\neq\phi_2(\phi(A))$. This proves \eqref{eq:twohydratwo}.

  Let $M$ be a 3-connected matroid having two inequivalent representations over $\GF(5)$. Then there exists a $\GF(5)\otimes\GF(5)$-matrix $A$ such that $M = M[I|A]$ and $\phi_1(A) \not \sim \phi_2(A)$.

  The restriction $\phi|_{\fun(\gauss)}:\fun(\gauss)\rightarrow\fun(\GF(5)\otimes\GF(5))$ is a bijection. If we apply Theorem~\ref{thm:lift} with lifting function $(\phi|_{\fun(\gauss)})^{-1}$ then Case~\ref{thm:lift}\eqref{enum:propfour} holds only for $\GF(5)\otimes\GF(5)$-matrices $A$ having a minor
 \begin{align}
   \begin{bmatrix}
     1 & 1 & 1\\
     1 & p & q
   \end{bmatrix} \textrm{ or }    \begin{bmatrix}
     1 & 1\\
     1 & p\\
     1 & q
   \end{bmatrix},
 \end{align}
 where $p,q \in \{(2,2), (3,3), (4,4)\}$. But Theorem~\ref{Utwofivestabilizer} implies that if $A$ has such a minor, then $\phi_1(A)$ and $\phi_2(A)$ will be strongly equivalent. Since both matrices have the same row and column indices, this implies $\phi_1(A) \sim \phi_2(A)$, a contradiction. Now \eqref{eq:onehydratwo} follows.
\end{proof}

\begin{theorem}\label{thm:gauss}
Let $M$ be a 3-connected matroid with a $U_{2,5}$- or $U_{3,5}$-minor. The following are equivalent:
\begin{enumerate}
  \item \label{eq:gaussfivefive}$M$ has 2 inequivalent representations over $\GF(5)$;
  \item \label{eq:gaussTU} $M$ is $\gauss$-representable;
  \item \label{eq:gaussevery}$M$ has two inequivalent representations over $\GF(5)$ and is representable over $\GF(p^2)$ for all primes $p \geq 3$ and over $\GF(p)$ when $p \equiv 1 \mod 4$.
\end{enumerate}
\end{theorem}
\begin{proof}
  \eqref{eq:gaussfivefive}$\Leftrightarrow$\eqref{eq:gaussTU} follows from the previous lemma. For \eqref{eq:gaussTU}$\Rightarrow$\eqref{eq:gaussevery} we use an argument similar to the proof of Theorem~\ref{thm:sru} where, as in the proof of Theorem~\ref{thm:splittable}, every element of $\gauss$ is of the form $2^k x$ for some $k \in \Z$, $x \in \Z[i]$.
  Finally, \eqref{eq:gaussfivefive} is a special case of \eqref{eq:gaussevery}.
\end{proof}

Let $\alpha$ be an indeterminate. For $k \geq 1$, a matroid is called \emph{$k$-cyclotomic} if it is representable over the partial field
\begin{align}
  \cyclo_k := \parf(\Q(\alpha), \{\alpha, \alpha-1, \alpha^2-1, \ldots, \alpha^k-1\}).
\end{align}
\begin{lemma}
  If $M$ is $\cyclo_k$-representable, then it is representable over every field that has an element $x$ whose multiplicative order is at least $k+1$. In particular, $M$ is representable over $\GF(q)$ for $q \geq k+2$.
\end{lemma}
\begin{proof}
  It is straightforward to construct a partial field homomorphism such that $\phi(\alpha) = x$.
\end{proof}
Let $\Phi_0(\alpha) := \alpha$ and let $\Phi_j$ be the $j$th \emph{cyclotomic polynomial}, i.e. the polynomial whose roots are exactly the primitive $j$th roots of unity. A straightforward observation is the following:
\begin{lemma}
$\cyclo_k = \parf(\Q(\alpha), \{\Phi_j(\alpha)\mid j = 0, \ldots, k\})$.
\end{lemma}
In particular $\cyclo_2 = \parf(\Q(\alpha), \{\alpha, \alpha-1,\alpha+1\})$.
\begin{lemma}$\fun(\cyclo_2) = \assoc\{1,\alpha,-\alpha,\alpha^2\}$. \end{lemma}
\begin{proof}
  Suppose $p := (-1)^s \alpha^x (\alpha-1)^y (\alpha^2-1)^z$ is a fundamental element.
  Every homomorphism $\phi:\cyclo_2\rightarrow \golrat$ and every homomorphism $\phi:\cyclo_2\rightarrow\gauss$ gives bounds on $x,y,z$. After combining several of these bounds a finite number of possibilities remains. We leave out the details.
\end{proof}
We conclude this section with the following result:
\begin{theorem}\label{thm:GFfourKtwo}
  Let $M$ be a matroid. The following are equivalent:
  \begin{itemize}
    \item $M$ is representable over $\GF(4)\otimes \gauss$;
    \item $M$ is representable over $\cyclo_2$.
  \end{itemize}
\end{theorem}
The proof consists, once more, of an application of Corollary~\ref{cor:pfequivalence}.

\section{An algebraic construction}\label{sec:rings}
With a theorem as general as the Lift Theorem, an interesting question becomes whether we can construct suitable partial fields $\hat \parf$ to which a given class of matroids lifts. In this section, we find the ``most general'' or ``algebraically most free'' partial field to which all $\parf$-representable matroids lift, a notion that we will make precise soon. Our starting point is Theorem~\ref{thm:pfinring}, which we prove now. For convenience we repeat the theorem here.
\begin{theorem}[Vertigan]\label{thm:pfinringrepeat}If $\parf$ is a partial field, then there exist a ring $\ring$ and a set $S\subseteq \ring^*$ such that $\parf \cong \parf(\ring,S)$.
\end{theorem}

\begin{proof}
  Let $\parf = (P, \oplus, \cdot, 0, 1_{\parf})$, and define $\group := (P\setminus \{0\}, \cdot, 1_{\parf})$. Recall that the group ring of $\group$ over $\Z$ is defined as
  \begin{align*}
     \Z[\group] := \{ \sum_{p \in \group} a_p \cdot p \mid a_p \in \Z, \textrm{finitely many } a_p \textrm{ are nonzero} \},
  \end{align*}
  where addition of two elements is componentwise and multiplication is defined by
  \begin{align}\label{eq:groupringdist}
    (\sum_{p\in\group} a_p \cdot p)(\sum_{p\in\group} b_p \cdot p) = \sum_{p,q\in \group} a_p b_q \cdot p q.
  \end{align}
  We identify $z \in \Z$ with $\sum_{i=1}^z 1_\parf$. We drop the $\cdot$ from the notation from now on. For clarity we write $p\oplus q$ if we mean addition in $\parf$, and $p+q$ if we mean (formal) addition in $\Z[\group]$. Consider the following subset of $\Z[\group]$:
  \begin{align*}
    V_1 := \{p + q \mid p \oplus q \doteq 0 \},
  \end{align*}
  and define the ideal $I_1 := V_1\Z[\group]$.
  \begin{claim}\label{cl:expandsum1}
    If $x \in I_1$ then $x = \pm s_1 \pm \cdots \pm s_k$ for some $s_1, \ldots, s_k \in V_1$.
  \end{claim}
  \begin{subproof}
    By definition $x = r_1 s_1 + \cdots + r_k s_k$ for $r_1, \ldots, r_k \in \Z[\group]$ and $s_1, \ldots, s_k \in V_1$. We consider one term.
    \begin{align*}
      r_i s_i = (\sum_{t \in \group} a_t t)(p+q) = \sum_{t \in \group} (a_t t(p+q)) = \sum_{t\in \group} (a_t (tp+tq)),
    \end{align*}
    where the last equality follows from \eqref{eq:groupringdist}. Since $p\oplus q \doteq 0$, also $tp \oplus tq \doteq 0$, by \eqref{ax:distributivity}. Hence $tp+tq \in V_1$. If $a_t > 0$ then
    \begin{align*}
      r_i s_i = \underbrace{(tp+tq) + \cdots + (tp+tq)}_{a_t \textrm{ terms}}.
    \end{align*}
    If $a_t < 0$ then
    \begin{align*}
      r_i s_i = \underbrace{-(tp+tq) - \cdots - (tp+tq)}_{-a_t \textrm{ terms}}.
    \end{align*}
    Summing over $i$ now yields the claim.
  \end{subproof}
  \begin{claim}
    $1_\parf \not \in I_1$.
  \end{claim}
  \begin{subproof}
    Suppose $1_\parf \in I_1$. By Claim~\ref{cl:expandsum1}, $1_\parf = \pm s_1 \pm \cdots \pm s_k$ for some $s_1, \ldots, s_k \in V_1$. We focus on the $s_i$ in which the coefficient of $1_\parf$ is not equal to 0. The only element of $V_1$ for which this holds is $1_\parf + (-1_\parf)$. It follows that, in $\pm s_1 \pm \cdots \pm s_k$, the coefficient of $(-1_\parf)$ is equal to that of $1_\parf$, which contradicts the assumption that $\pm s_1 \pm \cdots \pm s_k = 1_\parf$.
  \end{subproof}
  Now let $\ring_1 := \Z[\group]/I_1$. Consider the following subset of $\ring_1$:
  \begin{align*}
    V_2 := \{p+q+r + I_1 \mid (p\oplus q) \oplus r \doteq 0\},
  \end{align*}
  and define the ideal $I_2 := V_2 \ring_1$.
  \begin{claim}\label{cl:expandsum2}
    If $x \in I_2$ then $x = s_1 + \cdots + s_k$ for some $s_1, \ldots, s_k \in V_2$.
  \end{claim}
  \begin{subproof}
    By definition $x = r_1 s_1 + \cdots + r_k s_k$ for $r_1, \ldots, r_k \in \ring_1$ and $s_1, \ldots, s_k \in V_2$. We consider one term.
    \begin{align*}
      r_i s_i & = (\sum_{t \in \group} a_t t)(p+q+u) + I_1  = \sum_{t \in \group} (a_t t(p+q+u)) + I_1\\
       & = \sum_{t\in \group} (a_t (tp+tq+tu)) + I_1.
    \end{align*}
    Since $(p\oplus q)\oplus u \doteq 0$, also $(tp \oplus tq)\oplus tu \doteq 0$, by \eqref{ax:distributivity}. Hence $tp+tq+tu+I_1 \in V_2$. If $a_t > 0$ then
    \begin{align*}
      r_i s_i = \underbrace{(tp+tq+tu) + \cdots + (tp+tq+tu)}_{a_r \textrm{ terms}} + I_1.
    \end{align*}
    If $a_t < 0$ then we observe that $-p + I_1 = (-p) + I_1$, and obtain
    \begin{align*}
      r_i s_i = \underbrace{((-tp)+(-tq)+(-tu)) + \cdots + ((-tp)+(-tq)+(-tu))}_{-a_t \textrm{ terms}} + I_1.
    \end{align*}
    Summing over $i$ now yields the claim.
  \end{subproof}
  Now let $\ring_2 := \ring_1/I_2$, $\group_2 := \langle \{p + I_1 + I_2 \mid p \in \group\}\rangle$, and define $\parf' := \parf(\ring_2, \group_2)$. Our aim is to prove $\parf \cong \parf'$. To that end we construct a partial field isomorphism. Let $\phi:\parf\rightarrow\parf'$ be defined by
  \begin{align*}
    \phi(p) := p + I_1 + I_2.
  \end{align*}
  \begin{claim}
    $\phi$ is a partial field homomorphism.
  \end{claim}
  \begin{subproof}
    For $p, q\in P$, $\phi(p)\phi(q) = (p+I_1+I_2)(q+I_1+I_2) = pq + I_1 + I_2 = \phi(pq)$. If $p,q,r \in P$ are such that $p\oplus q \doteq r$ then $\phi(p) + \phi(q) = p + q + I_1 + I_2 = -(-r) + I_1 + I_2 = r + I_1 + I_2 = \phi(p\oplus q)$, since $p+q+(-r) \in V_2$ and $r + (-r) \in V_1$. Clearly $r + I_1 + I_2 \in \group_2 \cup \{0\}$, so $\phi(p) + \phi(q) \doteq \phi(r)$.
  \end{subproof}
  \begin{claim}\label{cl:apfpfbijection}
    $\phi$ is a bijection.
  \end{claim}
  \begin{subproof}
    Obviously $\phi$ is surjective. Suppose $p,q \in P$ are such that $p \neq q$ yet $\phi(p) = \phi(q)$. Then $p-q + I_1 \in I_2$. By Claim~\ref{cl:expandsum2}, $p-q = s_1 + \cdots + s_k$ for some $s_1, \ldots, s_k \in V_2$. For each $s_i$, pick representatives $p_i, q_i, r_i \in P$ such that $s_i = p_i + q_i + r_i + I_1$ and $(p_i\oplus q_i)\oplus r_i \doteq 0$. Define the multiset
    \begin{align*}
      S := \bigcup_{i=1}^k \{p_i,q_i,r_i\}.
    \end{align*}
    We build two associations for $S$. First, since $(p_i \oplus q_i) \oplus r_i \doteq 0$ and $0\oplus 0 \doteq 0$, we can build an association whose root node is labelled by $0$. Second, pick an $s \in S$. The only elements of $S$ contributing to the coefficient of $s + I_1$ in $s_1 + \cdots + s_k$ are $s$ and $(-s)$. Hence, for each $s \in S \setminus \{p,(-q)\}$, there is an element $(-s) \in S\setminus \{p,q\}$. By repeatedly pairing these elements we can build a pre-association where the children of the root node are labelled $p$ and $(-q)$. But the associative law then implies $p\oplus (-q) \doteq 0$, and hence $p = q$, contradicting our assumption.
  \end{subproof}
  In particular, Claim~\ref{cl:apfpfbijection} implies that $\phi$ is nontrivial.
  \begin{claim}
    $\phi$ is an isomorphism.
  \end{claim}
  \begin{subproof}
    Let $p,q,r \in P$ be such that $p+q+I_1+I_2 = r + I_1 + I_2$. We have to show that $p\oplus q \doteq r$.
    Since $p+q + (-r) + I_1 \in I_2$, there are $s_1, \ldots, s_n \in V_2$ such that $p+q + (-r) + I_1 = s_1 + \cdots + s_n$. For each $s_i$, pick representatives $p_i, q_i, r_i \in P$ such that $s_i = p_i + q_i + r_i + I_1$ and $(p_i\oplus q_i)\oplus r_i \doteq 0$. Define the multiset
    \begin{align*}
      S := \{r\}\cup \bigcup_{i=1}^k \{p_i,q_i,r_i\}.
    \end{align*}
    Using the same argument as in the previous claim we construct two pre-associations for $S$: one where the children of the root node are $r, 0$, and one where the children of the root node are $p, q$. Since $r\oplus 0 \doteq r$, the result follows from the associative law.
  \end{subproof}
  With this claim the proof is complete.
\end{proof}

Note that we have proven that $\parf \cong \parf(\ring_2, \group_2)$, not $\parf \cong \parf(\ring_2,\ring_2^*)$. It could be that $\group_2$ is a strict subgroup of $\ring_2^*$. 

\begin{corollary}\label{cor:pfimpliesfield}
  If $M$ is representable over a partial field $\parf$ then $M$ is representable over a field.
\end{corollary}
\begin{proof}
  Let $\parf = \parf(\ring,S)$, and let $A$ be a $\parf$-matrix such that $M = M[I|A]$. If every $x \in \ring\setminus 0$ is invertible then $\ring$ is a field. If some $x \in \ring\setminus 0$ is not invertible then $x\ring$ is a proper ideal of $\ring$. A standard result from commutative ring theory implies the existence of a maximal ideal $I \supseteq x\ring$, and then $\ring/I$ is a field (see, for example, Page 2 of Matsumura~\cite{Mat86}). There is a nontrivial ring homomorphism $\phi:\ring\rightarrow\ring/I$, and therefore, by Corollary~\ref{cor:hom}, $M = M[I|\phi(A)]$.
\end{proof}
Clearly every ring homomorphism yields a partial field homomorphism. On the  other hand, not all partial field homomorphisms extend to ring homomorphisms. The following example shows this. Let $\ring := \GF(2)\times\GF(7)$, and let $\parf := \GF(2)\otimes\GF(7)$. Let $\phi: \parf \rightarrow \reg$ be determined by $\phi(1,1)=\phi(1,2)=\phi(1,4) = 1$ and $\phi(1,6)=\phi(1,5)=\phi(1,3)=-1$. This is a partial field homomorphism. However, in $\ring$ we have $(1,2)+(1,4) = (1,3)+(1,3) = (0,6)$. It follows that $\phi$ cannot be extended to a homomorphism $\phi':\ring\rightarrow\Q$. The following theorem overcomes this problem. Recall from Definition~\ref{def:subpf} that $\parf[S]$ is the sub-partial field of $\parf$ with multiplicative group generated by $-1$ and $S$.

\begin{theorem}Let $\parf$, $\parf'$ be partial fields such that $\parf = \parf[\fun(\parf)]$ and $\parf' = \parf[\fun(\parf')]$, and suppose $\phi:\parf\rightarrow\parf'$ is a partial field homomorphism. Then there exist rings $\ring$, $\ring'$ and sets $S\subseteq \ring^*$, $S' \subseteq (\ring')^*$, such that $\parf \cong \parf(\ring,S)$, $\parf'\cong \parf(\ring',S')$, and such that $\phi$ can be extended to a ring homomorphism $\phi':\ring\rightarrow\ring'$. 
\end{theorem}
\begin{proof}
  Let $\ring_2$, $\ring_2'$ be the rings constructed in the proof of Theorem~\ref{thm:pfinringrepeat}. Every element of $\parf$ can be expressed as a product of fundamental elements and $-1$. From this it follows that there exists a ring homomorphism $\phi'': \Z[\parf_1^*] \rightarrow \ring_2'$. But $I_1+I_2 \subseteq \ker(\phi'')$. It follows that there exists a well-defined homomorphism $\phi':\ring_2\rightarrow\ring_2'$.
\end{proof}
The restriction on $\parf$, $\parf'$ in this theorem is rather light, as the following propositions show. We prove the first in \cite{PZ08conf}. The main idea is to look at induced cycles in the bipartite graph of a normalized representation.
\begin{proposition}\label{cor:closureispf}
  If a matroid $M$ is representable over a partial field $\parf$, then $M$ is representable over $\parf[\fun(\parf)]$.
\end{proposition}
\begin{proposition}\label{lem:homparfsfun}
  Let $\parf_1$, $\parf_2$ be partial fields and $\phi:\parf_1\rightarrow\parf_2$ a partial field homomorphism. Then there exists a partial field homomorphism $\phi': \parf_1[\fun(\parf_1)]\rightarrow\parf_2[\fun(\parf_2)]$.
\end{proposition}
\begin{proof}
  Let $\parf_1' := \parf_1[\fun(\parf_1)]$ and let $\parf_2' := \parf_2[\fun(\parf_2)]$. Then $\phi' := \phi|_{\parf_1'}: \parf_1'\rightarrow\parf_2$ is a partial field homomorphism. Clearly $\phi(-1) = -1$. Let $p = p_1\cdots p_k \in \parf_1'$, where $p_1, \ldots, p_k \in \fun(\parf_1')$. Then $\phi(p) = \phi(p_1)\cdots\phi(p_k) \in \parf_2'$. Hence the image of $\phi'$ is contained in $\parf_2'$, which completes the proof.
\end{proof}

Now that we can embed a partial field in a ring, we are ready for a construction of partial fields $\hat\parf$ satisfying the conditions of Corollary~\ref{cor:pfequivalence}.
\begin{definition}\label{def:levelpf}
  Let $\parf$ be a partial field.
  We define the \emph{lift} of $\parf$ as
  \begin{align}
     \lift\parf := \parf(\ring_{\parf}/I_{\parf}, \tilde F_{\parf}),
  \end{align}
  where $\tilde F_{\parf} := \{\tilde p \mid p \in \fun(\parf)\}$ is a set of indeterminates, one for every fundamental element, $\ring_{\parf} := \Z[\tilde F]$ is the polynomial ring over $\Z$ with indeterminates $\tilde F_{\parf}$, and $I_{\parf}$ is the ideal generated by the following polynomials in $\ring_{\parf}$:
  \begin{enumerate}
    \item\label{lev:zeroone} $\tilde 0 - 0$; $\tilde 1 - 1$;
    \item\label{lev:minone} $\tilde{-1} + 1$ if $-1 \in \fun(\parf)$;
    \item\label{lev:add} $\tilde p+\tilde q-1$, where $p,q \in \fun(\parf)$, $p + q \doteq 1$;
    \item\label{lev:mul} $\tilde p\tilde q-1$, where $p,q \in \fun(\parf)$, $pq = 1$;
    \item\label{lev:threeterm} $\tilde p\tilde q\tilde r-1$, where $p,q,r \in \fun(\parf)$, $pqr = 1$.
  \end{enumerate}
\end{definition}
We show that a matroid is $\parf$-representable if and only if it is $\lift\parf$-representable. First we need a lemma.
\begin{lemma}\label{lem:lifthom}
  Let $\parf$ be a partial field. There exists a nontrivial partial field homomorphism $\phi:\lift\parf\rightarrow\parf$ such that $\phi(\tilde p + I_{\parf}) = p$ for all $p \in \fun(\parf)$.
\end{lemma}
\begin{proof}
  Let $\ring$ be a ring such that $\parf = \parf(\ring,S)$ for some $S$. Then $\psi:\ring_{\parf}\rightarrow\ring$ determined by $\psi(\tilde p) = p$ for all $\tilde p \in \tilde F_{\parf}$ is obviously a ring homomorphism. Clearly $I_{\parf} \subseteq \ker(\psi)$, so $\phi':\ring_{\parf}/I_{\parf}\rightarrow \ring$ determined by $\phi'(\tilde p + I_{\parf}) = \psi(p)$ for all $\tilde p \in \tilde F_{\parf}$ is a well-defined ring homomorphism. Then $\phi := \phi'|_{\lift\parf}$ is the desired partial field homomorphism. Since $1 \not \in I_{\parf}$, $\phi$ is nontrivial.
\end{proof}
\begin{lemma}\label{lem:liftpf}
  Let $\parf$ be a partial field. A matroid is $\parf$-representable if and only if it is $\lift\parf$-representable.
\end{lemma}
\begin{proof}
  Let $\hat\parf := \lift\parf$ and let $\phi$ be the homomorphism from Lemma~\ref{lem:lifthom}. We define $^\uparrow: \fun(\parf)\rightarrow \fun(\hat\parf)$ by $p^\uparrow = \tilde p + I_{\parf}$. By \ref{def:levelpf}\eqref{lev:add},\eqref{lev:mul} this is a lifting function for $\phi$. Now all conditions of Corollary~\ref{cor:pfequivalence} are satisfied.
\end{proof}
The partial field $\lift\parf$ is the most general partial field for which the lift theorem holds, in the following sense:
\begin{theorem}
  Suppose $\parf$, $\hat\parf$, $\phi$, $^\uparrow$ are such that all conditions of Corollary~\ref{cor:pfequivalence} are satisfied. Then there exists a nontrivial homomorphism $\psi:\lift\parf\rightarrow\hat\parf$.
\end{theorem}
\begin{proof}
  Let $\psi':\ring_\parf\rightarrow\hat\parf$ be determined by $\psi'(\tilde p) = p^\uparrow$ for all $p \in \fun(\parf)$. This is clearly a ring homomorphism. But since all conditions of Corollary~\ref{cor:pfequivalence} hold, $I_{\parf} \subseteq \ker(\psi')$. It follows that there exists a well-defined homomorphism $\psi:\lift\parf\rightarrow\hat\parf$ as desired.
\end{proof}
Homomorphisms between lifts of partial fields are more well-behaved than homomorphisms between arbitrary partial fields:
\begin{lemma}
Let $\parf_1$, $\parf_2$ be partial fields, and let $\ring_{\parf_1}/I_{\parf_1}$, $\ring_{\parf_2}/I_{\parf_2}$ be the rings as in Definition~\ref{def:levelpf}. Let $\phi_i:\lift\parf_i \rightarrow\parf_i$ be the homomorphisms from Lemma~\ref{lem:lifthom}. Suppose that there exists a nontrivial partial field homomorphism $\phi:\parf_1 \rightarrow \parf_2$. Then there exists a nontrivial partial field homomorphism $\psi:\lift\parf_1\rightarrow\lift\parf_2$ that is the restriction of a ring homomorphism $\ring_{\parf_1}/I_{\parf_1} \rightarrow\ring_{\parf_2}/I_{\parf_2}$, such that the following diagram commutes:
\begin{align}
\begin{CD}
\lift\parf_1 @>\psi>> \lift\parf_2\\
@V\phi_1VV @VV\phi_2V\\
\parf_1 @>\phi>> \parf_2
\end{CD}
\end{align}
\end{lemma}
\begin{proof}
  We define $\psi': \ring_{\parf_1}\rightarrow\ring_{\parf_2}/I_{\parf_2}$ by $\psi'(\tilde p) = \tilde q + I_{\parf_2}$, where $\tilde q$ is such that $\phi(p) = q$. Again, this is obviously a ring homomorphism, and $I_{\parf_1} \subseteq \ker(\psi')$. The homomorphism $\psi:\ring_{\parf_1}/I_{\parf_1}\rightarrow \ring_{\parf_2}/I_{\parf_2}$ determined by $\psi(\tilde p +I_{\parf_1}) = \psi'(\tilde p)$ is therefore well-defined. The diagram now commutes by definition, and therefore nontriviality of $\psi$ follows from that of $\phi$.
\end{proof}

The importance of Lemma~\ref{lem:liftpf} is that we can now \emph{construct} partial fields for which the conditions of Corollary~\ref{cor:pfequivalence} hold. We use algebraic tools such as Gr\"obner basis computations over rings to get insight in the structure of $\lift\parf$. In particular, we adapted the method described by Baines and V\'amos~\cite{BV03} to verify the claims in Table~\ref{tab:levelpf}.

\begin{table}[tbp]
  \begin{center}
  \begin{tabular}{llll}
    \toprule
    $\parf\phantom{X}$ & $\GF(2)\otimes\GF(3)$ & $\GF(3)\otimes\GF(4)$ & $\GF(3)\otimes\GF(5)$  \\
    \addlinespace
    $\lift\parf$ & $\reg$ & $\psru$ & $\dyadic$  \\
    \midrule
    $\parf$ & $\GF(3)\otimes\GF(7)$ & $\GF(3)\otimes\GF(8)$ & $\GF(4)\otimes\GF(5)$  \\
    \addlinespace
    $\lift\parf$ & $\splittable$ & $\nreg$ & $\golrat$  \\
    \midrule
    $\parf$ & $\GF(5)\otimes\GF(7)$ & $\GF(5)\otimes\GF(8)$ & $\GF(4)\otimes\GF(5)\otimes\GF(7)$\\
    \addlinespace
    $\lift\parf$ & $\GF(5)\otimes\GF(7)$ & $\GF(5)\otimes\GF(8)$ & $\golrat\otimes\GF(7)$\\
    \bottomrule
  \end{tabular}
  \end{center}
  \caption{Some lifts of partial fields.}\label{tab:levelpf}
\end{table}

The obvious question is now: is $\lift\parf \not \cong \parf$ for other choices of $\parf = \GF(q_1)\otimes\cdots\otimes\GF(q_k)$? The last three entries in Table~\ref{tab:levelpf} indicate that sometimes the answer is negative. In these finite fields there seem to be relations that enforce $\lift\parf \cong \parf$. But Theorems~\ref{thm:gauss} and \ref{thm:GFfourKtwo} indicate that there are other uses still for the Lift Theorem. We conclude this section with a modification of Definition~\ref{def:levelpf} that accommodates the characterization of the Gaussian partial field.

\begin{definition}
  Let $\parf$ be a partial field and ${\cal A}$ a set of $\parf$-matrices.
  We define the ${\cal A}$-\emph{lift} of $\parf$ as
  \begin{align}
     \lift_{\cal A}\parf := \parf(\ring_{\parf}/I_{\parf}, \tilde F_{\parf}),
  \end{align}
  where $\tilde F_{\parf} := \{\tilde p \mid p \in \fun(\parf)\}$ is a set of symbols, one for every fundamental element, $\ring_{\parf} := \Z[\tilde F]$ is the polynomial ring over $\Z$ in indeterminates $\tilde F_{\parf}$, and $I_{\parf}$ is the ideal generated by the following polynomials in $\ring_{\parf}$:
  \begin{enumerate}
    \item $\tilde 0 - 0$; $\tilde 1 - 1$;
    \item $\tilde{-1} + 1$ if $-1 \in \fun(\parf)$;
    \item $\tilde p+\tilde q-1$, where $p,q \in \fun(\parf)$, $p + q \doteq 1$;
    \item $\tilde p\tilde q-1$, where $p,q \in \fun(\parf)$, $pq = 1$;
    \item $\tilde p\tilde q\tilde r-1$, where $p,q,r \in \fun(\parf)$, $pqr = 1$, and
    \begin{align}
      \begin{bmatrix}1 & 1 & 1\\ 1 & p & q^{-1}
      \end{bmatrix} \minorof A
    \end{align}
    for some $A \in {\cal A}$.
  \end{enumerate}
\end{definition}
We omit the proof of the following lemma.
\begin{lemma}
  Let $\parf$ be a partial field and ${\cal A}$ a set of $\parf$-matrices, and let $M$ be a matroid. If $M = M[I|A]$ for some $A \in {\cal A}$ then $M$ is $\lift_{\cal A}\parf$-representable.
\end{lemma}

\section{A number of questions and conjectures}\label{sec:questions}
While writing this paper we asked ourselves numerous questions. To some the answer can be found in this paper or in \cite{PZ08conf}, but in this section we present a few that are still open.

Theorems such as those in Section~\ref{sec:examples} show the equivalence between representability over infinitely many fields and over a finite number of finite fields. The following conjecture generalizes the characterization of the near-regular matroids:
\begin{conjecture}\label{con:finfieldlist}
  Let $k$ be a prime power. There exists a number $n_k$ such that, for all matroids $M$, $M$ is representable over all fields with at least $k$ elements if and only if it is representable over all finite fields $\GF(q)$ with $k \leq q \leq n_k$.
\end{conjecture}
To our disappointment the techniques in the present paper failed to prove this conjecture even for $k=4$. We offer the following candidate:
\begin{conjecture}\label{con:cyclokall}
  A matroid $M$ is representable over all finite fields with at least $4$ elements if and only if $M$ is representable over
  \begin{align}
    \parf_4 := \parf(\Q(\alpha),\{\alpha,\alpha-1,\alpha+1,\alpha-2\}),
  \end{align}
  where $\alpha$ is an indeterminate.
\end{conjecture}
Originally we posed this conjecture with $\cyclo_2$ instead of $\parf_4$. This would imply that all such matroids have at least two inequivalent representations over $\GF(5)$. But consider $M_{8591} := M[I|A_{8591}]$, where $A_{8591}$ is the following $\parf_4$-matrix:
\begin{align}
  A_{8591} := \begin{bmatrix}
    1 & 1 & 0 & \alpha & 1\\
    0 & 1 & 1 & \alpha & \alpha^{-1}\\
    1 & 0 & \alpha & \alpha & 1\\
    0 & 0 & 1 & 1 & 0
  \end{bmatrix}.
\end{align}
This matroid was found by Royle in Mayhew and Royle's catalog of small matroids~\cite{MR08} as a matroid representable over $\GF(4)$, $\GF(7)$, $\GF(8)$, and uniquely representable over $\GF(5)$. Hence $M_{8591}$ is not representable over $\cyclo_2$ (a fact that can be proven using tools from our forthcoming paper~\cite{PZ08conf}).

\begin{question}
  To what extent is a partial field $\parf$ determined by the set of finite fields $\GF(q)$ for which there exists a homomorphism $\phi:\parf \rightarrow \GF(q)$?
\end{question}
The previous example shows that $\parf$ is certainly not uniquely determined: both $\cyclo_2$ and $\parf_4$ have homomorphisms to all finite fields with at least $4$ elements, but $M_{8591}$ is only representable over the latter.

\begin{question}
  Are there systematic methods to determine the full set of fundamental elements for (certain types of) partial fields?
\end{question}
Semple~\cite{Sem97} determined the set of fundamental elements for a class of partial fields that he calls the $k$-regular partial fields. In this paper we computed $\fun(\parf)$ using ad hoc techniques, the only recurring argument being the fact that a homomorphism $\phi:\parf\rightarrow\parf'$ maps $\fun(\parf)$ to $\fun(\parf')$. We give two further illustrations.
First, consider the partial field
\begin{align}
  \ger := \parf(\Q,\{2,3\}).
\end{align}
This innocent-looking partial field, an extension of the dyadic partial field, has a finite number of fundamental elements, the least obvious of which are obtained from the relations $2^2-3=1$ and $3^2-2^3=1$. That there is indeed no other such relation is a classical but nonobvious result. It was proven by Gersonides in 1342
(see, for example, Peterson~\cite{Pet99} for a modern exposition). Consideration of $\parf(\Q, \{x,y\})$ for other pairs $x,y$ brings us into the realm of Catalan's Conjecture. This conjecture was posed more than 150 years ago and settled only in 2002.

Second, consider the partial field
\begin{align}
  \uniform_1^{(2)} := \parf(\GF(2)(\alpha), \{\alpha, 1+\alpha\}).
\end{align}
$\fun(\uniform_1^{(2)})$ has infinite size, since $\alpha^{2^k} - 1 = (\alpha + 1)^{2^k}$ for all $k \geq 0$.

The partial field $\lift\parf$ gives information about the representability of the set of $\parf$-representable matroids over other fields. An interesting question is how much information it gives.
\begin{question}\label{con:levelnatural}
  Which partial fields $\parf$ are such that whenever the set of $\parf$-representable matroids is also representable over a field $\field$, there exists a homomorphism $\phi:\lift\parf\rightarrow\field$?
\end{question}
In \cite{PZ08conf} we will show that each of $\reg$, $\psru$, $\dyadic$, $\nreg$, $\splittable$, $\golrat$, $\gauss$ has this property.
\begin{question}
  Let $\phi:\lift\parf\rightarrow\parf$ be the canonical homomorphism. For which partial fields $\parf$ is $\phi|_{\fun(\lift\parf)}:\fun(\lift\parf)\rightarrow\fun(\parf)$ a bijection?
\end{question}
This bijection exists for all examples in this paper and results in an obvious choice of lifting function. If there is always such a bijection then it is not necessary to introduce an abstract lifting function. In that case the proof of the Lift Theorem can be simplified to some extent. A related conjecture is the following:
\begin{conjecture}
  $\lift^2\parf \cong \lift\parf$.
\end{conjecture}

We end with a conjecture that seems to be only just outside the scope of the Lift Theorem:
\begin{conjecture}
  A matroid is representable over $\GF(2^k)$ for all $k > 1$ if and only if it is representable over $\uniform_1^{(2)}$.
\end{conjecture}

In an earlier version of this paper we also conjectured that a matroid is representable over $\GF(4)\otimes\R$ if and only if it is representable over $\golrat$. Afterwards we found that the Pappus matroid is a counterexample to this.

\paragraph*{Acknowledgements}We thank Hendrik Lenstra for suggesting the $k$-Cyclotomic partial field. We also thank Christian Eggermont for some helpful comments on rings of integers in algebraic number fields. We thank Gordon Royle for his quick and friendly responses when we asked him for data from the catalog of small matroids~\cite{MR08}. His examples prevented the authors from embarking on several wild goose chases. Finally we thank two anonymous referees for carefully reading the paper, and for making useful suggestions that improved the quality of the final version.

\appendix

\section{When should we call a sum ``defined''?}\label{sec:counterex}
The notion of a sum $p_1 + \cdots + p_n$ being \emph{defined} appears somewhat complicated. Semple and Whittle~\cite{SW96} give a simpler definition: $p_1 + \cdots + p_n$ is defined if there exists some association of $\{p_1, \ldots, p_n\}$. Unfortunately, this simpler definition has a problem. Consider the following matrices:
\begin{align}
  A :=
   \begin{bmatrix}
  1 & 1 & 0 & 0 & 0\\
  1 & 0 & 1 & 0 & 0\\
  1 & 0 & 0 & 1 & 0\\
  0 & b + a & c & d - a & -1\\
  0 & -a & 0 & a & 1
  \end{bmatrix},
   & \quad  B :=
     \begin{bmatrix}
  1 & 1 & 0 & 0 & 0\\
  1 & 0 & 1 & 0 & 0\\
  1 & 0 & 0 & 1 & 0\\
  0 & b & c & d & 0\\
  0 & -a & 0 & a & 1
  \end{bmatrix},
\end{align}
where $B$ is obtained from $A$ by adding the last row to the next to last.
Then $\det(A) = (b + a) + c + (d - a) - a + a$ and $\det(B) = b + c + d$. In both sums no cancellation has taken place: all terms missing from the formal determinant are 0. Now consider the following instantiation over $\ring := \Z/51\Z$:
\begin{align}
    a = 37, b = 7, c = 23, d = 11.
\end{align}
Then none of $b + c$, $b+d$, $c+d$ are invertible, yet $a$, $b$, $c$, $d$, $1$, $-1$, $(b+a)$, $((b+a)+c)$, $d-a$, $((b+a)+c) + (d-a)$ are. It follows that in $\parf(\ring,\ring^*)$, $\det(A)$ is defined in the sense of Semple and Whittle~\cite{SW96}, whereas $\det(B)$ is not.

This is a counterexample to Proposition \ref{prop:matops}\eqref{eq:matopsfour}, which is therefore false under the old definition. This proposition is used for pretty much everything that comes after it in Semple and Whittle~\cite{SW96}, so it is important to find a way to fix it. The proposed change in the meaning of a sum being \emph{defined} is one way to do that. To make absolutely sure that this is indeed the case, we give a proof of Proposition \ref{prop:matops} using the new definition.
\begin{proof}[Proof of Proposition \ref{prop:matops}]
Assume $B$ was obtained from $A$ by transposition. Then
\begin{align}
    \det(B) & = \sum_{\sigma \in S_n} \sign(\sigma) b_{1\sigma(1)} b_{2\sigma(2)}\cdots b_{n\sigma(n)}\\
            & = \sum_{\sigma \in S_n} \sign(\sigma) a_{\sigma(1)1} a_{\sigma(2)2}\cdots a_{\sigma(n)n}
\end{align}
which is nothing but a permutation of the terms of $\det(A)$.

Assume $B$ was obtained from $A$ by swapping rows $1$ and $2$. Then
\begin{align}
    \det(B) & = \sum_{\sigma \in S_n} \sign(\sigma) b_{1\sigma(1)} b_{2\sigma(2)} b_{3\sigma(3)}\cdots b_{n\sigma(n)}\\
            & = \sum_{\sigma \in S_n} \sign(\sigma) a_{2\sigma(1)} a_{1\sigma(2)} a_{3\sigma(3)}\cdots a_{n\sigma(n)}\\
            & = \sum_{\sigma' \in S_n} \sign(\sigma') a_{2\sigma'(2)} a_{1\sigma'(1)} a_{3\sigma'(3)}\cdots a_{n\sigma'(n)}
\end{align}
where $\sigma' = \sigma \circ (1,2)$ (in cycle notation; cycles act from the right). Therefore $\sign(\sigma') = - \sign(\sigma)$, from which the second part of the proposition follows.

For the third part, assume we multiply row 1 by a constant $p$. Then
\begin{align}
    \det(B) & = \sum_{\sigma \in S_n} \sign(\sigma) b_{1\sigma(1)} b_{2\sigma(2)} \cdots b_{n\sigma(n)}\\
            & = \sum_{\sigma \in S_n} \sign(\sigma) p a_{1\sigma(1)} a_{2\sigma(2)} \cdots a_{n\sigma(n)}\\
            & = p \det(A).
\end{align}
Here the last line follows from Axiom~\eqref{ax:distributivity}.

For the final part we prove the following, more general lemma:
\begin{lemma}Let $Z := \{1,\ldots, n\}$. Let $A = [a | X]$ and $B = [b | X]$ be $Z\times Z$ matrices with entries in $\parf$ such that $A[Z, Z\setminus 1] = B[Z, Z\setminus 1] = X$. If $\det(A)$, $\det(B)$, $\det(A)+\det(B)$ and all entries of the vector $a+b$ are defined, then $\det([a+b | X]) \doteq \det(A) + \det(B)$.
\end{lemma}
\begin{subproof}
Set $C = [a+b|X]$. Then
\begin{align}
    &\det(C) & = & \sum_{\sigma \in S_n} \sign(\sigma) c_{1\sigma(1)} c_{2\sigma(2)}\cdots c_{n\sigma(n)}\\
    &        & = & \sum_{\sigma \in S_n} \sign(\sigma) (a+b)_{1\sigma(1)} c_{2\sigma(2)}\cdots c_{n\sigma(n)}\\
    &        & = & \sum_{\sigma \in S_n} \sign(\sigma) (a+b)_{1\sigma(1)} c_{2\sigma(2)}\cdots c_{n\sigma(n)}\notag\\
    &        &   & - \sum_{\sigma \in S_n} \sign(\sigma) b_{1\sigma(1)} c_{2\sigma(2)}\cdots c_{n\sigma(n)}\notag\\
    &        &   & + \sum_{\sigma \in S_n} \sign(\sigma) b_{1\sigma(1)} c_{2\sigma(2)}\cdots c_{n\sigma(n)}\\
    &        & = & \sum_{\sigma \in S_n} \sign(\sigma) a_{1\sigma(1)} c_{2\sigma(2)}\cdots c_{n\sigma(n)}\notag\\
    &        &   & + \sum_{\sigma \in S_n} \sign(\sigma) b_{1\sigma(1)} c_{2\sigma(2)}\cdots c_{n\sigma(n)}.\label{eq:detsumthing}
\end{align}
For \eqref{eq:detsumthing} we used the fact that, if $(a+b)$ is defined, then $(a+b)-b \doteq a$ (an easy consequence of Axioms~\eqref{ax:zero} and \eqref{ax:associativelaw}), together with Axiom~\eqref{ax:distributivity}. For the final expression it is easy to provide an association: take associations $T_A$, $T_B$ for $\det(A)$, $\det(B)$; add a new root vertex $r$ and edges $r_Ar$, $r_Br$. This is a pre-association for $\det(C)$. Since $r_A$ is labelled by $\det(A)$ and $r_B$ by $\det(B)$, we have that $r$ is labelled by $\det(A)+\det(B)$, which was defined by assumption.
\end{subproof}
Returning to the proof of the proposition, let $B$ be obtained from $A$ by adding row $i$ to row $1$, where we assume that $a_{1j}+a_{ij}$ is defined for all $j$. Let $A'$ be the matrix obtained by replacing the first row of $A$ by the $i$th row, and leaving all other rows unaltered. Since the first and the $i$th row of $A'$ are identical, $\det(A')=0$ (it is easy to find an association, since the terms of the determinant cancel pairwise). Applying the lemma to $A$, $A'$ we conclude that $\det(B) \doteq \det(A) + \det(A') = \det(A)$, as desired.
\end{proof}

Since the proposed change occurs at the fringes of the definitions related to partial fields, it does not cause much damage. In fact, all other propositions, lemmas and theorems of \cite[{Sections 1--6}]{SW96} are true under the new definition.

As a final remark we note that, even with our definition, the following occurs. Consider the sum $1+1+1$ in $\ring := \Z/4\Z$. The units of this ring are $1$, $3$, and the only nontrivial sum that is defined in $\parf(\ring,\ring^*)$ is $1+3\doteq 0$. It follows that $1+1+1$ is undefined in $\parf(\Z/4\Z,(\Z/4\Z)^*)$ yet a unit in $\ring$.

\section{A catalog of partial fields}
\begin{figure}[htbp]
  \center
  \includegraphics{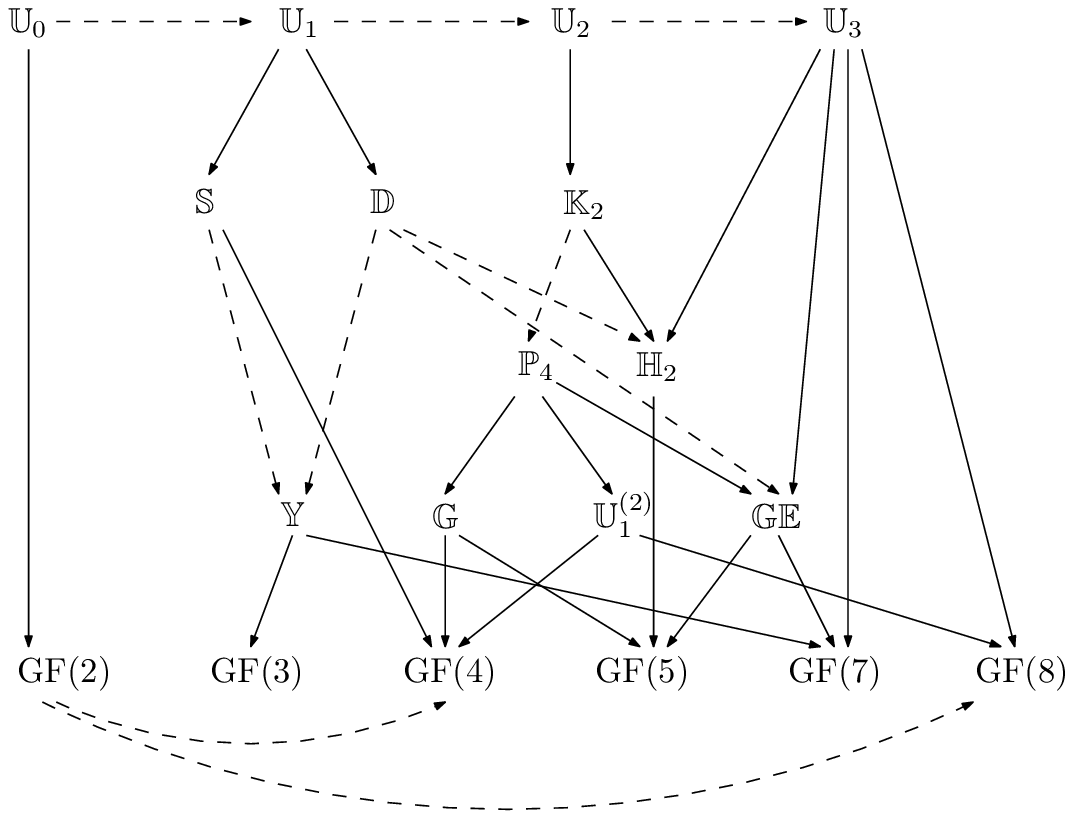}
  \caption{Some partial fields and their homomorphisms. A (dashed) arrow from $\parf'$ to $\parf$ indicates that there is an (injective) homomorphism $\parf'\rightarrow \parf$.}\label{fig:pfcat}
\end{figure}

In this appendix we summarize all partial fields introduced in this paper, as well as a class of partial fields introduced by Semple~\cite{Sem97}. Like rings, partial fields form a category. The regular partial field, $\reg$, has a homomorphism to every other partial field. In Figure~\ref{fig:pfcat} we display the relations between the partial fields from this appendix. Recall from Definition~\ref{prop:pffromring} that $\parf(\ring, S)$ is the partial field $(\langle S\cup \{-1\}\rangle, +, \cdot, 0, 1)$, where multiplication and addition are the restriction of the operations in $\ring$. $\fun(\parf)$ denotes the set of fundamental elements of $\parf$, and $\assoc F$ is as defined in Section~\ref{ssec:funds}.

\begin{description}
  \item[The regular partial field,] $\reg$:
  \begin{itemize}
    \item $\reg = \parf(\Q,\{-1,0,1\})$;
    \item $\fun(\reg) = \{0,1\}$;
    \item There is a homomorphism to every partial field $\parf$;
    \item Isomorphic to $\lift(\GF(2)\otimes\GF(3))$;
    \item There are finitely many excluded minors for $\reg$-representability (Theorem~\ref{thm:regex}).
  \end{itemize}
  \item[The near-regular partial field,] $\nreg$:
  \begin{itemize}
    \item $\nreg = \parf(\Q(\alpha), \{\alpha,1-\alpha\})$, where $\alpha$ is an indeterminate;
    \item $\fun(\nreg) = \assoc\{1,\alpha\} = \left\{0,1,\alpha,1-\alpha,\frac{1}{1-\alpha},\frac{\alpha}{\alpha-1},\frac{\alpha-1}{\alpha},\frac{1}{\alpha}\right\}$;
    \item There is a homomorphism to every field with at least three elements;
    \item Isomorphic to $\lift(\GF(3)\otimes\GF(8))$ and $\lift(\GF(3)\otimes\GF(4)\otimes\GF(5))$;
    \item There are finitely many excluded minors for $\nreg$-representability \cite{HMZ09}.
  \end{itemize}
  \item[The $k$-uniform partial field,] $\uniform_k$:
  \begin{itemize}
    \item $\uniform_k = \parf\left(\Q(\alpha_1, \ldots, \alpha_k), \left\{p-q \mid p,q \in \{0,1,\alpha_1,\ldots,\alpha_k\}, p \neq q\right\}\right)$, where $\alpha_1,\ldots,\alpha_k$ are  indeterminates;
    \item Introduced by Semple~\cite{Sem97} as the \emph{$k$-regular} partial field;
    \item Semple~\cite{Sem97} proved that
    \begin{align}
      \fun(\uniform_k) = & \left\{ \frac{a-b}{c-b} \,\Big|\, a,b,c \in \{0,1,\alpha_1,\ldots,\alpha_k\}, \textrm{ distinct}\right\} \cup\notag\\
      & \left\{ \frac{(a-b)(c-d)}{(c-b)(a-d)} \,\Big|\, a,b,c,d \in \{0,1,\alpha_1,\ldots,\alpha_k\}, \textrm{ distinct}\right\};
    \end{align}
    \item There is a homomorphism to every field with at least $k+2$ elements~\cite{Sem97};
    \item Finitely many excluded minors for $\uniform_k$-representability are $\uniform_{k'}$-representable for some $k' > k$~\cite{OSV00}.
  \end{itemize}
  \item[The sixth-roots-of-unity ($\sru$) partial field,] $\psru$:
  \begin{itemize}
    \item $\psru = \parf(\C, \zeta)$, where $\zeta$ is a root of $x^2-x+1=0$.
    \item $\fun(\psru) = \assoc\{1,\zeta\} = \{0,1,\zeta,1-\zeta\}$;
    \item There is a homomorphism to $\GF(3)$, to $\GF(p^2)$ for all primes $p$, and to $\GF(p)$ when $p \equiv 1 \mod 3$;
    \item Isomorphic to $\lift(\GF(3)\otimes\GF(4))$;
    \item There are finitely many excluded minors for $\psru$-representability \cite{GGK}.
  \end{itemize}
  \item[The dyadic partial field,] $\dyadic$:
  \begin{itemize}
    \item $\dyadic = \parf(\Q,2)$;
    \item $\fun(\dyadic) = \assoc\{1,2\} = \{0,1,-1,2,1/2\}$;
    \item There is a homomorphism to every field that does not have characteristic two;
    \item Isomorphic to $\lift(\GF(3)\otimes\GF(5))$.
  \end{itemize}
  \item[The union of $\sru$ and dyadic,] $\splittable$:
  \begin{itemize}
    \item $\splittable = \parf(\C,\{2,\zeta\})$, where $\zeta$ is a root of $x^2-x+1=0$;
    \item $\fun(\splittable) = \assoc\{1,2,\zeta\} = \{0,1,-1,2,1/2,\zeta,1-\zeta\}$;
    \item There is a homomorphism to $\GF(3)$, to $\GF(p^2)$ for all odd primes $p$, and to $\GF(p)$ when $p \equiv 1 \mod 3$;
    \item Isomorphic to $\lift(\GF(3)\otimes\GF(7))$.
  \end{itemize}
  \item[The $2$-cyclotomic partial field,] $\cyclo_2$:
  \begin{itemize}
    \item $\cyclo_2 = \parf(\Q(\alpha), \{\alpha, \alpha-1,\alpha+1\})$, where $\alpha$ is an indeterminate;
    \item $\fun(\cyclo_2) = \assoc\{1,\alpha,-\alpha,\alpha^2\}$;
    \item There is a homomorphism to $\GF(q)$ for $q \geq 4$;
    \item Isomorphic to $\lift(\GF(4)\otimes\gauss)$.
  \end{itemize}
  \item[The $k$-cyclotomic partial field,] $\cyclo_k$:
  \begin{itemize}
      \item $\cyclo_k = \parf(\Q(\alpha), \{\alpha, \alpha-1, \alpha^2-1, \ldots, \alpha^k-1\})$, where $\alpha$ is an indeterminate;
      \item $\cyclo_k = \parf(\Q(\alpha), \{\Phi_j(\alpha)\mid j = 0, \ldots, k\})$, where $\Phi_0(\alpha) = \alpha$ and $\Phi_j$ is the $j$th \emph{cyclotomic polynomial};
      \item There is a homomorphism to $\GF(q)$ for $q \geq k+2$.
  \end{itemize}
  \item[The Gersonides partial field,] $\ger$:
  \begin{itemize}
    \item $\ger = \parf(\Q, \{2,3\})$;
    \item $\fun(\ger) = \assoc\{1,2,3,4,9\}$;
    \item There is a homomorphism to every field that does not have characteristic two or three.
  \end{itemize}
  \item[The partial field ] $\parf_4$:
  \begin{itemize}
    \item $\parf_4 = \parf(\Q(\alpha), \{\alpha,\alpha-1,\alpha+1,\alpha-2\})$, where $\alpha$ is an indeterminate;
    \item $\fun(\parf_4) = \assoc\{1,\alpha,-\alpha,\alpha^2,\alpha-1,(\alpha-1)^2\}$;
    \item There is a homomorphism to every field with at least four elements.
  \end{itemize}
  \item[The Gaussian partial field,] $\gauss$:
  \begin{itemize}
    \item $\gauss = \parf(\C,\{i,1-i\})$, where $i$ is a root of $x^2+1=0$;
    \item $\fun(\gauss) = \assoc\{1,2,i\} = \left\{0,1,-1,2,\tfrac{1}{2},i,i+1,\tfrac{i+1}{2},1-i,\tfrac{1-i}{2},-i\right\}$;
    \item There is a homomorphism to $\GF(p^2)$ for all primes $p \geq 3$, and to $\GF(p)$ when $p \equiv 1 \mod 4$;
    \item A matroid is $\gauss$-representable if and only if it is dyadic or has at least two inequivalent $\GF(5)$-representations.
  \end{itemize}
  \item[The near-regular partial field modulo two,] $\nreg^{(2)}$:
  \begin{itemize}
    \item $\nreg^{(2)} = \parf(\GF(2)(\alpha), \{\alpha, 1+\alpha\})$, where $\alpha$ is an indeterminate;
    \item $\fun(\nreg^{(2)}) = \{0,1\}\cup \assoc \left\{\alpha^{2^k} \mid k \in \N \cup \{0\}\right\}$;
    \item There is a homomorphism to $\GF(2^k)$ for all $k \geq 2$.
  \end{itemize}
  \item[The golden ratio partial field,] $\golrat$:
  \begin{itemize}
    \item $\golrat = \parf(\R, \tau)$, where $\tau$ is the positive root of $x^2-x-1=0$;
    \item $\fun(\golrat) = \assoc\{1,\tau\} = \{0,1,\tau, -\tau, 1/\tau, -1/\tau, \tau^2, 1/\tau^2\}$;
    \item There is a homomorphism to $\GF(5)$, to $\GF(p^2)$ for all primes $p$, and to $\GF(p)$ when $p \equiv \pm 1 \mod 5$;
    \item Isomorphic to $\lift(\GF(4)\otimes\GF(5))$.
  \end{itemize}
\end{description}

\bibliography{matbib}
  \bibliographystyle{svzwambib}
\end{document}